\documentclass[11pt,reqno]{amsart}
\usepackage{amsmath,amssymb,amsthm,mathtools,mathrsfs}
\usepackage[T1]{fontenc}
\usepackage{lmodern,microtype}
\usepackage{booktabs,array,enumitem,xcolor}
\usepackage{needspace,etoolbox}
\usepackage[a4paper,left=25mm,right=25mm,top=26mm,bottom=27mm,
  headheight=14pt,headsep=8mm,footskip=12mm]{geometry}
\usepackage{xurl}
\usepackage{hyperref}
\hypersetup{colorlinks=true,linkcolor=black,citecolor=black,urlcolor=blue!45!black,pdftitle={The Woo--Yong conjecture},pdfauthor={Minghua Dou}}
\patchcmd{\abstract}{\normalfont\Small}{\normalfont\small}{}{%
  \PackageError{paper-layout}{Unable to set abstract size}{}}
\numberwithin{equation}{section}
\newtheoremstyle{paperplain}
  {8pt plus 2pt minus 1pt}{7pt plus 2pt minus 1pt}
  {\itshape}{}{\bfseries}{.}{.6em}{}
\newtheoremstyle{paperdefinition}
  {8pt plus 2pt minus 1pt}{7pt plus 2pt minus 1pt}
  {\normalfont}{}{\bfseries}{.}{.6em}{}
\newtheoremstyle{paperremark}
  {8pt plus 2pt minus 1pt}{7pt plus 2pt minus 1pt}
  {\normalfont}{}{\itshape}{.}{.6em}{}
\theoremstyle{paperplain}
\newtheorem{theorem}{Theorem}[section]
\newtheorem{lemma}[theorem]{Lemma}
\newtheorem{proposition}[theorem]{Proposition}
\newtheorem{corollary}[theorem]{Corollary}
\newtheorem{conjecture}[theorem]{Conjecture}
\theoremstyle{paperdefinition}

\newtheorem{example}[theorem]{Example}
\theoremstyle{paperremark}
\newtheorem{remark}[theorem]{Remark}
\newcommand{\C}{\mathbb C}
\newcommand{\Z}{\mathbb Z}
\newcommand{\Q}{\mathbb Q}
\newcommand{\A}{\mathbb A}
\newcommand{\PP}{\mathbb P}
\newcommand{\Gm}{\mathbb G_m}
\newcommand{\OO}{\mathcal O}
\newcommand{\M}{\mathcal M}
\DeclareMathOperator{\Cl}{Cl}
\DeclareMathOperator{\Pic}{Pic}
\DeclareMathOperator{\Spec}{Spec}
\DeclareMathOperator{\Sing}{Sing}
\DeclareMathOperator{\NGor}{NGor}
\DeclareMathOperator{\NonQGor}{NonQGor}
\DeclareMathOperator{\NonCart}{NonCart}
\DeclareMathOperator{\Tor}{Tor}
\DeclareMathOperator{\Dem}{Dem}
\DeclareMathOperator{\rank}{rank}
\DeclareMathOperator{\codim}{codim}
\DeclareMathOperator{\id}{id}
\DeclareMathOperator{\Gr}{Gr}
\DeclareMathOperator{\OG}{OG}
\DeclareMathOperator{\LG}{LG}
\DeclareMathOperator{\Frac}{Frac}
\DeclareMathOperator{\GL}{GL}
\DeclareMathOperator{\SL}{SL}
\DeclareMathOperator{\Irr}{Irr}

\DeclareMathOperator{\ind}{ind}

\DeclareMathOperator{\Hom}{Hom}
\DeclareMathOperator{\NC}{NC}
\newcommand{\Bbad}{B_{\mathrm{bad}}}
\newcommand{\tors}{\mathrm{tors}}
\newcommand{\sh}{\mathrm{sh}}
\newcommand{\eps}{\varepsilon}
\newcommand{\red}{\mathrm{red}}
\AtBeginDocument{%
  \setlength{\parindent}{1.3em}%
  \setlength{\parskip}{2.5pt plus .7pt minus .4pt}%
  \setlength{\abovedisplayskip}{9pt plus 2pt minus 2pt}%
  \setlength{\belowdisplayskip}{9pt plus 2pt minus 2pt}%
  \setlength{\abovedisplayshortskip}{5pt plus 2pt minus 1pt}%
  \setlength{\belowdisplayshortskip}{7pt plus 2pt minus 2pt}%
}
\makeatletter
\renewcommand\subsection{\@startsection{subsection}{2}%
  \z@{.9\baselineskip\@plus.3\baselineskip}{.3\baselineskip}%
  {\normalfont\bfseries}}
\makeatother
\BeforeBeginEnvironment{theorem}{\Needspace{5\baselineskip}}
\BeforeBeginEnvironment{lemma}{\Needspace{5\baselineskip}}
\BeforeBeginEnvironment{proposition}{\Needspace{5\baselineskip}}
\BeforeBeginEnvironment{corollary}{\Needspace{5\baselineskip}}
\BeforeBeginEnvironment{conjecture}{\Needspace{5\baselineskip}}
\BeforeBeginEnvironment{example}{\Needspace{5\baselineskip}}
\AtBeginEnvironment{thebibliography}{\raggedright}
\setlist[enumerate]{itemsep=4pt,topsep=6pt,parsep=2pt}
\renewcommand{\arraystretch}{1.12}
\allowdisplaybreaks[1]
\title[The Woo--Yong conjecture]{The Woo--Yong conjecture}
\author{Minghua Dou}
\address{The University of Manchester, Manchester, United Kingdom}
\email{minghua.dou@outlook.com}
\date{November 6, 2026}
\begin{document}
\begin{abstract}
We prove the Woo--Yong conjecture: the non-Gorenstein locus of a type~A Schubert variety is the union of its generically non-Gorenstein maximal singular components. More strongly, at every scheme point, local Weil divisor classes are detected at the generic points of the maximal singular components containing that point. For arbitrary finite crystallographic type, the analogous localization map is injective on the torsion subgroup of the source local class group. This determines Cartier indices of locally $\Q$-Cartier divisors and gives the Woo--Yong criterion on the $\Q$-Gorenstein locus. A type~D family shows the necessity of this hypothesis: the ambient variety is generically Gorenstein along its singular locus but not $\Q$-Gorenstein along a deeper nonempty locus. The type~A proof and the counterexamples are geometric; the general torsion proof also uses bounded finite inputs with explicit certificates and reproducible reconstruction programs.
\end{abstract}
\maketitle
\tableofcontents
\clearpage
\setcounter{section}{-1}
\section{A survey of the geometry behind the conjecture}\label{sec:survey}

This section explains the objects and the obstruction measured by the conjecture, with particular attention to the different scopes of the two detection arguments. Readers familiar with Schubert varieties and local divisor theory may proceed to Section~\ref{sec:introduction}. The proofs in the rest of the paper do not depend on the survey.

The motivating question concerns singular spaces with finitely many natural strata: does the ambient local geometry at a general point of each maximal singular component determine where a given algebraic defect occurs, including on smaller strata? Here the defect is failure to be Gorenstein. On a Schubert variety this means that the canonical sheaf is not a line bundle. Divisor class groups encode that failure, and localization supplies the maps by which it is tested.

\subsection{Flags, Schubert cells, and Bruhat order}
A complete flag in $\C^n$ is a chain $0=F_0\subset F_1\subset\cdots\subset F_n=\C^n$ with $\dim F_i=i$. These flags form the smooth projective variety $\GL_n/B$, where $B$ is the group of invertible upper triangular matrices. Fix the standard flag $E_j=\langle e_1,\ldots,e_j\rangle$. The $B$-orbits record relative positions with respect to $E_\bullet$ and are the Schubert cells \cite[Section~1.2]{Bri05}.

The cell $C_w$ indexed by $w\in S_n$ contains the coordinate flag $F_i=\langle e_{w(1)},\ldots,e_{w(i)}\rangle$. Row elimination identifies it with $\A^{\ell(w)}$, where $\ell(w)$ counts inversions. Its closure $X_w=\overline{C_w}$ is a Schubert variety. Explicitly, it is defined by
\[
 \dim(F_i\cap E_j)\ge\#\{a\le i:w(a)\le j\}.
\]
In matrix charts these incidence inequalities become rank bounds, and hence equations given by minors \cite[Section~3]{WY08}. The cells are smooth, while their closures may be singular.

Bruhat order is defined by $u\le w$ if $C_u\subseteq X_w$, equivalently $X_u\subseteq X_w$; thus $X_w=\coprod_{u\le w}C_u$. A cover $u\lessdot w$ has $\ell(w)=\ell(u)+1$. In $S_2$, the identity labels a point and the other permutation labels the open affine cell in $\PP^1$. In general the longest element labels the whole flag variety. Bruhat order is a partial order: a dimension inequality alone does not imply a closure inclusion \cite[Chapter~2]{BB05}.

A partial flag remembers only specified subspaces; a Grassmannian remembers one. Forgetting steps gives $G/B\to G/P$, a locally trivial projection with smooth flag variety fibre. This allows local singularity questions to pass between complete and partial flags; Section~4 gives the compatible detector maps. Symplectic and orthogonal groups give flags compatible with bilinear forms. Their Weyl groups replace ordinary permutations. Types~B and~C have the same abstract signed-permutation Weyl group but different root lengths, which produce different integral character matrices below.

\subsection{What a local ring and a generic point remember}
For $X=\Spec R$, a point is a prime ideal $\mathfrak p$ and its local ring is $R_{\mathfrak p}$. At a closed complex point this ring consists of rational functions defined on some neighbourhood. General prime ideals are equally important \cite[Chapter~II, Sections~2--3]{Har77}.

An irreducible closed subvariety $Z$ has a generic point $\eta_Z$ with closure $Z$. The ambient ring $\OO_{X,\eta_Z}$ records transverse behaviour along a nonempty open subset of $Z$: its dimension is $\codim_X Z$, while its residue field is $\C(Z)$.

\Needspace{5\baselineskip}
For example, for the line $Z=V(y)$ in $\A^2$, one has
\[
 \OO_{\A^2,\eta_Z}=\C[x,y]_{(y)}\simeq\C(x)[y]_{(y)},
\]
a discrete valuation ring with parameter $y$ and residue field $\C(x)$. In contrast, the intrinsic ring $\OO_{Z,\eta_Z}$ is the field $\C(x)$. Every generic-point test here concerns the ambient ring.

If $x\in Z$, then $\eta_Z$ is a generalization of $x$ and there is a localization map $\OO_{X,x}\to\OO_{X,\eta_Z}$. Freeness of a module persists under this localization. Its converse generally fails: a property may hold at general points of a component and fail at special ones. The conjecture asks for a specific converse for a particular obstruction on Schubert varieties.

This also explains the field extensions used in Section~4. A generic point along a positive-dimensional cell has a function field as residue field, even though the variety is defined over $\C$. The local arguments therefore retain an arbitrary extension $F\supset\C$ and use the same $F$ on every closed subslice participating in one argument.

\subsection{Regularity, normality, and Gorensteinness}
A Noetherian local ring $(R,\mathfrak m)$ is regular if $\dim_{R/\mathfrak m}(\mathfrak m/\mathfrak m^2)=\dim R$. For varieties over $\C$, regularity is the local algebraic form of smoothness. A singular point may nevertheless satisfy weaker conditions that control duality and divisors.

An integral variety is normal if its local rings are integrally closed domains. In codimension one these rings are discrete valuation rings, so rational functions have integral orders of vanishing along prime divisors. Normality also permits functions and reflexive sheaves to be recovered across complements of codimension at least two \cite[Chapter~II, Section~6]{Har77}\cite[Tags 031S, 0EBJ]{The26}.

A local ring is Cohen--Macaulay if its depth equals its dimension. Depth is the maximal length of a regular sequence in the maximal ideal, each element being a nonzerodivisor after the preceding ones have been set equal to zero. Complete intersections are Cohen--Macaulay, but the latter class is much larger \cite[Chapters~1--2]{BH93}.

The canonical module expresses local duality. On a smooth $d$-dimensional variety, the canonical sheaf is the line bundle $\bigwedge^d\Omega_X^1$. On a normal variety it is its rank-one reflexive extension from the smooth locus. For the Cohen--Macaulay local rings considered here, Gorensteinness is equivalent to the canonical module being free of rank one \cite[Chapter~3]{BH93}\cite[Tag 0DW9]{The26}. Thus a singular ring can retain duality with one local generator. The standard implications \cite[Tag 0AWV]{The26} are
\[
 \text{regular}\ \Longrightarrow\ \text{local complete intersection}
 \ \Longrightarrow\ \text{Gorenstein}
 \ \Longrightarrow\ \text{Cohen--Macaulay}.
\]
Normality is a separate condition. Schubert varieties in this paper are normal and Cohen--Macaulay \cite{Ram85,Bri05}. Our question is where their canonical sheaves stop being line bundles, not merely where the varieties become singular.

\subsection{Divisors and the group containing the obstruction}
A prime Weil divisor on a normal integral variety is an irreducible closed subvariety of codimension one. Weil divisors are finite integer combinations of these; a rational function $f\ne0$ defines $\operatorname{div}(f)=\sum_E\operatorname{ord}_E(f)E$. Divisors modulo principal divisors form $\Cl(X)$. A Cartier divisor is locally the divisor of one nonzero rational function, with unit ratios on overlaps, and determines a line bundle. Every Weil divisor on a smooth variety is Cartier, whereas this can fail on a normal singular variety \cite[Chapter~II, Section~6]{Har77}\cite{Fos73}.

The local group $\Cl(\OO_{X,x})$ measures the failure near $x$: $D$ is Cartier at $x$ precisely when its local class is zero. The qualification ``local'' matters. A point on $\PP^1$ defines the globally nontrivial line bundle $\OO_{\PP^1}(1)$, although its class is zero in every local class group. Globally, $\Cl(\PP^1)\simeq\Z$.

Equivalently, a Weil class is a rank-one reflexive sheaf $\mathcal F$, with $\mathcal F\simeq\mathcal F^{**}$ and $\mathcal F^*=\mathcal Hom(\mathcal F,\OO_X)$. The group operation is tensor product followed by double dual \cite[Tags 0EBJ, 0EBM]{The26}. Cartierness at a point is freeness of rank one over the local ring. This description supplies a closed nonfree locus; the flat descent used below is proved in Lemma~\ref{alg:flat}.

Let $K_X$ represent the canonical sheaf. Since $X$ is Cohen--Macaulay, it is Gorenstein at $x$ if and only if $[K_X]=0$ in $\Cl(\OO_{X,x})$.
A normal local ring is factorial precisely when its class group is zero. Factoriality therefore tests every class; Gorensteinness tests the distinguished canonical class. A divisor is $\Q$-Cartier at $x$ if some positive multiple is Cartier, meaning that its local class is torsion. Its Cartier index is the order of that class. The analogous condition for $K_X$ is $\Q$-Gorensteinness. Nonzero finite-order and infinite-order canonical classes both obstruct Gorensteinness, but a torsion-detection theorem applies only to the former.

\subsection{Three model calculations}
\begin{example}[A singular Gorenstein point with infinite class group]\label{ex:node}
The rank-at-most-one $2\times2$ matrix variety $D_{2,2}=\Spec\C[a,b,c,d]/(ad-bc)$ has dimension three. The gradient $(d,-c,-b,a)$ vanishes only at its vertex, where the tangent space has dimension four. It is the ordinary threefold node; being a hypersurface, it is Gorenstein.

The ruling primes $P=(a,b)$ and $Q=(a,c)$ satisfy $\operatorname{div}(a)=P+Q$. Inverting $a$ gives the factorial ring $\C[a^{\pm1},b,c]$, whose units are $\C^*a^{\Z}$. Divisor localization therefore identifies the global class group with
\[
 (\Z[P]\oplus\Z[Q])/\Z(P+Q)\simeq\Z.
\]
Its generator stays nonzero in the vertex local ring: if a homogeneous symbolic power became principal there, graded Nakayama would make it principal globally. Section~3 proves this comparison in detail; see also \cite[Section~8]{BV88}. Consequently the local class of $a_0P+b_0Q$ is $(a_0-b_0)[P]$. A singular Gorenstein point can carry many non-Cartier Weil divisors, and comparing the two ruling coefficients is an integral test.
\end{example}

\begin{example}[Finite-order divisors and the canonical class]\label{ex:veronese}
The invariant ring for simultaneous negation on $\A^2$ is
\[
 R=\C[s^2,st,t^2]\simeq\C[u,v,w]/(uw-v^2).
\]
It is normal and Gorenstein. At its vertex the prime $P=(u,v)$ is not principal, while $\operatorname{div}(u)=2P$ and $R[u^{-1}]=\C[u^{\pm1},v]$. Thus the vertex local class group is $\Z/2$: there is a torsion divisor obstruction, although the canonical class is zero.

For simultaneous negation on $\A^3$, the invariant ring is $\C[x_ix_j:1\le i\le j\le3]$. Its vertex local class group is $\Z/2$ \cite[Section~1.2, Proposition~I]{Wat74}. The involution acts on $dx_1\wedge dx_2\wedge dx_3$ by $-1$; the group has no pseudo-reflections and is not contained in $\SL_3$, so Watanabe's local criterion gives non-Gorensteinness \cite[Theorem~2]{Wat74}. The canonical class corresponds to the determinant character \cite[Section~2, Theorem~$1^\prime$, pp.~383--384]{Wat74}. It is the nonzero class, and its square is trivial, giving canonical index two. Equivalently, the toric lattice is the even-sum sublattice of $\Z^3$; the canonical divisor has odd total parity. These scalar Veronese cones are the type~C models studied in Section~5.
\end{example}

\begin{example}[An infinite-order canonical obstruction]\label{ex:determinantal}
For $r,s\ge2$, the rank-at-most-one variety $D_{r,s}$ of $r\times s$ matrices is normal and Cohen--Macaulay, has dimension $r+s-1$, and has vertex local class group $\Z$. Its row and column rulings represent primitive opposite classes. Up to the choice of sign for the generator, its canonical class is $(r-s)$ times that generator \cite[Corollary~8.4, Theorem~8.8 and Corollary~8.9]{BV88}. Indeed, the variety is the cone over $\PP^{r-1}\times\PP^{s-1}$ embedded by $\OO(1,1)$; the canonical bundle $\OO(-r,-s)$ maps to the difference of its coefficients after quotienting by the diagonal hyperplane class. The homogeneous local comparison is the same as in Example~\ref{ex:node}. Thus $D_{2,3}$ has a nonzero canonical class of infinite order and is not $\Q$-Gorenstein at its vertex.
\end{example}

These models have isolated singularities. Their punctured spectra are smooth, so Weil classes are line bundles on the puncture. For a nonisolated singularity a reflexive module need not be invertible on the whole puncture. The torsion argument first passes to the generic point of a component of the nonfree locus; only then does it obtain the necessary line bundle.

\subsection{From global classification to a locus problem}
A permutation contains a pattern $\pi\in S_m$ if some $m$ entries have the same relative order as $\pi$. Lakshmibai--Sandhya proved that a type~A Schubert variety is smooth precisely when its permutation avoids $3412$ and $4231$ \cite{LS90}\cite[Theorem~7.6]{WY23}. This converts a geometric property into a finite combinatorial test, but it does not determine Gorensteinness: Example~\ref{ex:node} is singular and Gorenstein.

The singular locus is closed and $B$-stable, hence a finite union of Schubert subvarieties. Its irreducible components are called maximal singular components, with maximality understood by inclusion, not dimension. Their generic singularities in type~A are described in \cite{BW03,Cor03,KLR03,Man01b,Man01a}. They supply a manageable list of ambient local models. The remaining issue is whether those models detect all failures on deeper strata.

Woo--Yong's global classification determines which complete type~A Schubert varieties are Gorenstein everywhere \cite{WY06}. Their Corollary~1 says that Gorensteinness at every maximal singular generic point implies Gorensteinness everywhere. The locus conjecture asks more: every non-Gorenstein point should lie on a maximal singular component that is already generically non-Gorenstein. With $\M(X)=\Irr(\Sing X)$, the asserted identity is
\begin{equation}\label{eq:survey-wy}
 \NGor(X)=\bigcup_{\substack{Z\in\M(X)\\\OO_{X,\eta_Z}\text{ not Gorenstein}}}Z.
\end{equation}
Openness gives the right-to-left inclusion. The other inclusion is substantive, particularly when generically Gorenstein and non-Gorenstein components coexist. Removing the known bad components does not permit another application of a theorem about complete Schubert varieties.

Subsequent work includes the interval-pattern approach \cite{WY08}, Perrin's determination of the Gorenstein locus in minuscule Schubert varieties \cite[Theorem~0.3 and Corollary~0.4]{Per09}, boundary blow-ups and a Demazure-subword criterion \cite{DS19}, and the survey \cite{WY23}. Global Cartier-class results \cite{LRY25} and canonical trace calculations for Grassmannian Schubert cycles and determinantal rings \cite{Kim26} address related questions. Section~\ref{sec:introduction} specifies their scopes.

\subsection{Why the conjecture leads to divisor detection}
At $x\in X$, localizing at the maximal singular generic points whose closures contain $x$ gives
\begin{equation}\label{eq:survey-rho}
 \rho_x:\Cl(\OO_{X,x})\longrightarrow
 \bigoplus_{\substack{Z\in\M(X)\\x\in Z}}\Cl(\OO_{X,\eta_Z}).
\end{equation}
The containment is necessary for these localization maps to exist. The canonical form of the conjecture asks whether $\rho_x([K_X])=0$ implies $[K_X]=0$. Injectivity of the whole map is stronger. Theorem~\ref{thm:mainA} proves that stronger assertion in type~A, including infinite-order classes, and consequently also determines the nonfactorial locus \cite[Conjecture~6.14]{WY08}\cite[Conjecture~7.39]{WY23}.

Theorem~\ref{thm:maintorsion} for arbitrary finite type asserts injectivity only on the source torsion subgroup. This qualification cannot be dropped: even the abstract map $\Z\to0$ is injective on its zero torsion subgroup while killing every nonzero integer. Likewise, vanishing of generic images does not make a local canonical class finite-order.

For a class already known to have finite order, injection preserves its order. The order in a finite direct sum is the least common multiple of the component orders, so a divisor known to be $\Q$-Cartier at $x$ satisfies the index formula of Corollary~\ref{cor:indices}. An empty least common multiple is one. The source hypothesis supplies finite order; the generic tests determine that finite order exactly.

\subsection{Transverse slices remove smooth directions}
Let $B^-$ be the opposite Borel and $C^u=B^-uB/B$. The Kazhdan--Lusztig slice $N_{u,w}=X_w\cap C^u$ has distinguished point $u$ and dimension $\ell(w)-\ell(u)$. There is an actual open product
\[
 X_w\cap\Omega_u\simeq C_u\times N_{u,w}
 \simeq\A^{\ell(u)}\times N_{u,w},
\]
compatible with Schubert strata \cite[Lemma~3.2]{WY08}\cite[Lemma~2.1]{KWY13}. The affine factor records movement along the cell, and the slice records transverse local geometry. Freeness of the relevant reflexive modules can be compared through the faithfully flat product.

In type~A, free matrix entries have root characters. Bruhat rank inequalities determine which root directions survive in the tangent space. A contracting positive grading then shows that the surviving entries generate the actual coordinate ring. This is why matrix algebra and permutation inequalities enter the proof together.

A tangent space records only first-order information, and a tangent cone is an associated graded object. Neither determines the original singularity: for instance, $y^2=x^3$ and $y^2=x^5$ have the same tangent cone $y^2=0$. Here the product and the determinantal models are actual scheme isomorphisms, carrying the prime divisors and their multiplicities. An unmarked generic-singularity classification alone would not identify which two divisors become the matrix rulings. Theorem~\ref{ta:filtered} supplies that placement information.

\subsection{The type A proof: comparing coefficients}
The node calculation illustrates the mechanism: when a matrix entry has two reduced ruling components, a divisor supported on them is locally trivial precisely when its two coefficients agree.

Choose the first-difference coordinate $t$ at $v<w$. Its nonvanishing locus is $N_{v,w}[t^{-1}]\simeq\Gm\times N_{u,w}$ with $\ell(u)=\ell(v)+1$. Induction handles divisor classes on this lower-dimensional slice.

A class vanishing on the open set can then be written $[D]=\sum_i c_i[E_i]$, where the $E_i$ are the prime components of $t=0$ and the $c_i$ are integers. This representation need not be unique. Reducedness gives the principal relation $\sum_i[E_i]=0$, so equality of all the coefficients suffices.

Some components come from Schubert boundary divisors and are called vertical. A possible additional component meets the dense cell and is called horizontal. For neighbouring active vertical components, a maximal element of a filtered Bruhat intersection gives a rank-one matrix cone in which the two chosen divisors are the rulings of one entry. The cone is singular only at its vertex, so its label corresponds to a maximal singular component of the ambient variety, an allowed detector.

The horizontal component requires an additional construction. Proposition~\ref{ta:horizontal} moves to a larger stratum while retaining a non-Cartier vertical divisor. A global divisor-class relation identifies the horizontal contribution in another two-ruling comparison. Thus every required neighbouring difference, including the last horizontal comparison, is detected. If every detector vanishes, all $c_i$ coincide and the principal relation gives $[D]=0$.

The proof therefore needs integral coefficients, reduced multiplicities, actual divisor markings, and globally maximal singular components. These properties are established separately; they are not inferred from a tangent-cone resemblance.

\subsection{The general-type proof: finite covers and integral degrees}
Finite-order classes admit an additional interpretation. If a line bundle $L$ on $U$ has a chosen trivialization $L^{\otimes m}\simeq\OO_U$, the powers of $L$ form a cyclic finite algebra. Locally its cover has equation $z^m=a$ with $a$ a unit, so in characteristic zero it is finite étale. Kummer theory relates these covers to torsion line bundles, with unit ambiguity $\Gamma(U,\OO_U)^*/\Gamma(U,\OO_U)^{*m}$ \cite{Gro71}\cite[Tag 03PK]{The26}.

To apply this, one first takes the generic point of an irreducible component of the nonfree locus. The reflexive module is then free on the punctured local spectrum. Passing to a strict henselization, a henselian local ring with separably closed residue field, allows simple roots to lift \cite[Tag 04GE]{The26}. Finite étale covers of the entire strict local spectrum split, although covers of its puncture may not. In the normal strictly henselian local rings of dimension at least two used here, units on the puncture extend and have $m$th roots. The unit ambiguity disappears. The descent argument uses injectivity for the relevant punctured Picard group; it does not identify the entire ordinary and strict local class groups.

The first detection method restricts covers to the Schubert boundary. Its top pieces correspond to coatoms $v\lessdot w$, and smaller Bruhat labels describe their overlaps. For $d=\ell(w)-\ell(u)$, the order complex of $(u,w)$ is a sphere $S^{d-2}$ \cite{BW82}. When $d\ge4$, its simple connectivity makes labels for split sheets on the pieces compatible. A local Lefschetz theorem then lifts the splitting from the boundary to the puncture \cite[Tag 0BLJ]{The26}. The restriction functor is used only for full faithfulness.

The second method uses a regular proper birational source, obtained by restricting a Bott--Samelson resolution to the slice. This restricted source is not itself asserted to be a tower of projective lines. Let $E_j$ be all of its prime divisors lying over the vertex. For a Kummer equation $z^m=a$, valuations off that fibre are divisible by $m$, so
\[
 \operatorname{div}(a)=mD+\sum_jc_jE_j.
\]
Taking degrees on actual proper integral cycles $Z_i$ in the fibre gives
\[
 N(c_j)_j\equiv0\pmod m,
 \qquad N_{ij}=\deg\bigl(\OO(E_j)|_{Z_i}\bigr).
\]
Smith normal form gives an annihilator for the torsion class. In particular, a full-column primitive matrix eliminates it: an integer left inverse exists, or equivalently the full-column maximal minors have gcd one. Rational rank is insufficient; the matrix $(2)$ permits a nonzero solution modulo $2$.

The geometry identifies every divisor column and realizes the character rows as integer combinations of genuine curves, using the Bott--Samelson fibre geometry of \cite{Ha04}. The effective Levi reduction bounds the required rank by slice dimension. In type~C the high-dimensional induction identifies a slice with nonzero strict puncture torsion as an actual scalar Veronese cone $\A^d/\{\pm1\}$. Its singularity is isolated, so it cannot hide torsion at a nonmaximal singular label. The bounded type~$F_4$ calculation and the type~$G_2$ word argument complete the finite-type result. The finite computations establish only their stated bounded inputs; geometry supplies the all-ranks and all-fields conclusions.

\subsection{The type D family and infinite-order obstruction}
For $r\ge k\ge3$, Section~\ref{sec:type-d} constructs a Schubert variety $X\subset\OG(k,2r+4)$ with a unique maximal singular component $S$ of codimension three. The transverse model at $\eta_S$ is the Gorenstein node. Nevertheless the ambient variety is non-$\Q$-Gorenstein at every point of a smaller locus $T\simeq\Gr(k,r)\subset S$, of codimension $2k$ in $X$. The Grassmannian $T$ is itself smooth; it is the ambient local rings of $X$ along $T$ that fail the canonical condition.

An incidence variety $I$ is a normal Gorenstein determinantal hypersurface and maps properly and birationally to $X$. It is isomorphic to $X$ away from $T$ and has fibre $\PP^1$ over every point of $T$, with no exceptional prime divisor. Along a contracted fibre curve $C$, the canonical line bundle of $I$ has degree $2k-4>0$.

If $mK_X$ were Cartier near its image, its pullback would have degree zero on $C$. Since there are no exceptional prime divisors, comparing canonical Weil divisors identifies that pullback with $mK_I$, whose degree is $m(2k-4)$. This contradiction applies to every scheme point of $T$, and proves an infinite-order local canonical obstruction.

The smallest case $(r,k)=(3,3)$ has type~$D_5$, dimension six, and $T$ a point. The case $(5,4)$ has type~$D_7$, dimension twelve, and $T\simeq\PP^4$. These orthogonal Grassmannians are nonminuscule, so the examples do not conflict with Perrin's theorem. They explain the additional term $\NonQGor(X)$ in Corollary~\ref{cor:indices}.

\subsection{Reading routes and background}
The three main assertions have distinct scopes: all divisor classes in type~A, finite-order classes in arbitrary finite type, and an infinite-order canonical obstruction in type~D. For the original conjecture, Sections~3 and~6 give the coefficient argument and its consequences. For the local-cover and lattice argument, Sections~4 and~5 give the geometry and induction. Section~7 is independent of the finite certificates and can be read directly after the divisor background.

For further background, \cite{Bri05} treats flag varieties and their divisor theory, while \cite{WY23} develops the combinatorial perspective. Canonical modules are treated in \cite{BH93}, divisor class groups in \cite{Har77,Fos73}, and determinantal models in \cite{BV88}. Bruhat topology is covered in \cite{BB05,BW82}, finite covers in \cite{Gro71}, and the normal affine semigroup methods in \cite{CLS11}. The new arguments begin with the marked detectors and their assembly into the local statements.

\section{Introduction}\label{sec:introduction}

Let $G$ be a connected reductive group over $\C$, let $B$ be a Borel subgroup, and let $P\supseteq B$ be a parabolic subgroup. A Schubert variety $X\subseteq G/P$ is normal and Cohen--Macaulay \cite{Ram85,Bri05}. Its canonical sheaf is therefore a rank-one reflexive sheaf, and its Gorenstein locus is the locus on which that sheaf is invertible \cite[Chapter~3]{BH93}. We ask whether this local obstruction can be detected at the generic points of the maximal singular components.

Write $\M(X)=\Irr(\Sing X)$, and denote the generic point of $Z\in\M(X)$ by $\eta_Z$. Define
\begin{equation}\label{eq:bad}
 \Bbad(X)=\bigcup_{\substack{Z\in\M(X)\\ \OO_{X,\eta_Z}\text{ is not Gorenstein}}} Z.
\end{equation}
This definition concerns the ambient local ring of $X$, rather than the intrinsic singularities of $Z$. Openness of the Gorenstein locus gives $\Bbad(X)\subseteq\NGor(X)$.

Woo and Yong classified globally Gorenstein Schubert varieties in type~A \cite{WY06}. Their Corollary~1 states that a type~A Schubert variety is globally Gorenstein if and only if its local rings at the generic points of all maximal singular components are Gorenstein. Their Conjecture~1 asks for the following pointwise strengthening.

\begin{conjecture}[Woo--Yong]\label{conj:wy}
If $X$ is a Schubert variety in a type~A flag variety, then $\NGor(X)=\Bbad(X)$.
\end{conjecture}

All equalities of loci in this paper are set-theoretic and include nonclosed points; they do not assert equality of scheme structures. Our convention is $X_w=\overline{BwB/B}$ and $\dim X_w=\ell(w)$; some original references use the opposite indexing convention.

The distinction between the global classification and Conjecture~\ref{conj:wy} is substantial. Removing some singular components usually leaves an open subset that is not a Schubert variety, so the global classification cannot simply be reapplied. The issue includes varieties with both generically Gorenstein and generically non-Gorenstein maximal singular components. A proof must handle points on the former components that lie on none of the latter.

\subsection{Divisor detection in type A}
For every scheme point $x\in X$, localization gives
\begin{equation}\label{eq:rho}
 \rho_x:\Cl(\OO_{X,x})\longrightarrow
 \bigoplus_{\substack{Z\in\M(X)\\x\in Z}}\Cl(\OO_{X,\eta_Z}).
\end{equation}
The condition $x\in Z$ makes $\eta_Z$ a generalization of $x$, so each component map is defined. An empty direct sum is the zero group.

\begin{theorem}[Type A divisor detection]\label{thm:mainA}
Let $X$ be a Schubert variety in a type~A complete or partial flag variety. For every scheme point $x\in X$, the map $\rho_x$ is injective. Equivalently, every Weil divisor $D$ on $X$ satisfies
\begin{equation}\label{eq:noncartier}
 \{x:D\text{ is not Cartier at }x\}
 =\bigcup_{\substack{Z\in\M(X)\\ D\text{ is not Cartier at }\eta_Z}} Z.
\end{equation}
\end{theorem}

Taking $D=K_X$ proves Conjecture~\ref{conj:wy}. Allowing all Weil divisors also gives the geometric nonfactorial-locus assertion of \cite[Conjecture~6.14]{WY08} and \cite[Conjecture~7.39]{WY23}. These consequences use the integral statement, including classes of infinite order.

The proof in Section~\ref{sec:type_a} inverts a coordinate associated with the first difference between two permutations. The resulting open subset is an actual product of a smaller Kazhdan--Lusztig slice and $\Gm$, and its complementary divisor is reduced.

Maximal elements of a filtered Bruhat intersection give rank-one determinantal slices in which two specified vertical components are the two rulings through one matrix entry. They detect successive integer differences of divisor coefficients. A separate prefix argument connects a horizontal component to the same comparison. Nagata localization and induction then prove injectivity. No cancellation of affine factors or identification of a tangent cone with the original slice is used.

\subsection{Finite-order detection in arbitrary type}
The same conclusion for all classes fails in other root systems. There is, however, a uniform assertion for finite-order classes.

\begin{theorem}[Torsion detection]\label{thm:maintorsion}
Let $X$ be a Schubert variety in $G/P$ for a finite crystallographic root system, possibly reducible. For every scheme point $x\in X$, the restriction of $\rho_x$ to $\Cl(\OO_{X,x})_{\tors}$ is injective.
\end{theorem}

The torsion hypothesis is imposed on the source class. A calculation after tensoring with $\Q$ would not establish this assertion.

\begin{corollary}[Indices and canonical loci]\label{cor:indices}
Suppose that $D$ is $\Q$-Cartier at $x$. Then
\[
 \ind_x(D)=\mathrm{lcm}_{Z\ni x}\ind_{\eta_Z}(D),
\]
where $Z$ ranges over $\M(X)$ and the empty least common multiple is $1$. On the $\Q$-Gorenstein locus, $X$ is Gorenstein at $x$ if and only if it is Gorenstein at every $\eta_Z$ with $x\in Z$. On all of $X$,
\[
 \NGor(X)=\NonQGor(X)\cup\Bbad(X).
\]
\end{corollary}

The proof of Theorem~\ref{thm:maintorsion} has a different geometric mechanism. In dimension at least four, finite covers of a strict local puncture are detected on its coatom boundary. The topology of an open Bruhat interval supplies compatible sheet labels; a local Lefschetz theorem lifts the resulting splitting. A second argument bounds torsion by degrees of actual integral cycles in a Bott--Samelson fibre.

Low-dimensional slices reduce to Levi subgroups of bounded rank. The simply laced and type~B cases follow from these arguments and low-rank bases. In type~C, the surviving strict torsion models are scalar second-Veronese cones, which have isolated singularities. Type~$F_4$ requires a finite certificate and type~$G_2$ an alternating-word argument. Throughout, field extensions of $\C$ are retained because a Schubert cell's generic point has a nonconstant residue field.

\subsection{The hypothesis beyond type A}
\begin{theorem}[A type D family]\label{thm:mainD}
For every $r\ge k\ge3$, there is a Schubert variety $X\subset\OG(k,2r+4)$ of dimension $k(r-k+2)$ whose singular locus is irreducible of codimension three and generically an ordinary threefold node. Nevertheless,
\[
 \NGor(X)=\NonQGor(X)\simeq\Gr(k,r),
 \qquad \codim_X\NGor(X)=2k.
\]
The same phenomenon occurs after passage to complete flags.
\end{theorem}

In this family $\Bbad(X)=\varnothing$. The smallest member has type~$D_5$, dimension six, and one non-Gorenstein point; the member $(r,k)=(5,4)$ has type~$D_7$, dimension twelve, and non-Gorenstein locus $\PP^4$. They contradict the unrestricted general-type formulation in \cite[Conjecture~0.1]{Per09}, while satisfying the scopes of the type~A and finite-order results: the obstructing canonical class has infinite order.

\subsection{Relation to earlier work}
Perrin proved the locus statement for minuscule Schubert varieties, including ordinary Grassmannian Schubert varieties \cite[Theorem~0.3 and Corollary~0.4]{Per09}. Woo--Yong's interval-pattern framework describes local properties through Bruhat intervals \cite{WY08}; their survey distinguishes the pointwise conjecture from the global classification and records verification through $S_6$ \cite[Conjecture~7.31 and Exercise~7.32]{WY23}.

The product charts used here come from \cite{WY08,KWY13}. Maximal singular components in type~A were classified independently in \cite{BW03,Cor03,KLR03,Man01b}; their generic local models are described in \cite{Cor03,Man01a}. Our filtered construction supplies the additional marked divisors required for coefficient comparisons.

Combining the local Grassmannian reduction for vexillary opposite Schubert varieties in \cite[Theorem~7.1 and Section~5]{AIJK23} with Perrin's theorem gives the corresponding covexillary locus result. Da Silva studies boundary blow-ups and gives a Demazure-subword criterion \cite[Theorem~0.2 and Corollary~0.7]{DS19}.

Li--Rietsch--Yang study global Cartier and $\Q$-Cartier classes on $G/P$ Schubert varieties \cite{LRY25}, addressing the global problem in \cite[Problem~8.24]{WY23}. Kimura determines canonical trace ideals and non-Gorenstein loci for Grassmannian Schubert cycles and determinantal rings \cite{Kim26}. These results have related but different scopes; none states the local injection \eqref{eq:rho} for all ordinary local rings of an arbitrary full-flag Schubert variety.

PermPAL \cite{PermPAL}, based on the framework of \cite{ABC26}, provides specifications and enumerations of permutation avoidance classes. Its class $\operatorname{Av}(3412,52341,635241)$ occurs in the iterated Grassmannian-bundle characterization of Alland--Richmond \cite{AR18}. This observation supplies combinatorial context and is not an input to either detection theorem.

The proof dependencies are explicit. In type~A, Theorem~\ref{ta:filtered} identifies the marked ruling divisors, Proposition~\ref{ta:horizontal} handles the horizontal component, and Lemma~\ref{ta:differences} and Theorem~\ref{ta:affine} complete the coefficient argument.

In arbitrary type, the Levi reduction, curve-degree construction, and boundary theorem provide the geometric inputs, while the type~C argument supplies the unbounded induction. The type~D construction is independent of these steps. Section~\ref{sec:reproducibility} and the accompanying source package specify the finite computations and their reconstruction.

\section{Local divisors and Schubert charts}\label{sec:algebra}

All schemes are Noetherian, and varieties are over $\C$ unless a field extension is specified. On a normal integral scheme, Weil divisor classes correspond to rank-one reflexive sheaves. A class vanishes in a local class group precisely when the corresponding module is free of rank one. We use the divisor--sheaf correspondence, the height-one criterion for reflexivity, and Serre's normality criterion \cite[Chapter~II, Section~6]{Har77}\cite[Tags 0EBJ, 0EBM, 0AVB, 031S]{The26}.

\begin{lemma}[Localization of divisor classes]\label{alg:localization}
Let $R$ be a normal Noetherian domain, let $0\ne t\in R$, and let $E_1,\ldots,E_s$ be the height-one primes containing $t$. Then
\[
 \bigoplus_{i=1}^s\Z[E_i]\longrightarrow\Cl(R)
 \longrightarrow\Cl(R[t^{-1}])\longrightarrow0
\]
is exact. Polynomial and Laurent polynomial extensions have the same class group as $R$.
\end{lemma}
\begin{proof}
Every prime divisor of $D(t)$ has a unique prime-divisor closure in $\Spec R$. If a divisor restricts to a principal divisor, subtract the divisor of the same rational function on $\Spec R$; the difference is supported on $V(t)$. This proves exactness, including the required kernel generation. For polynomial invariance, divisors dominating $\Spec R$ become principal over $\Frac(R)[T]$, while vertical valuations are extensions of the height-one valuations of $R$; Gauss's lemma identifies the relations. Localizing in $T$ gives Laurent invariance \cite[Chapter~II, Propositions~6.5--6.6]{Har77}\cite{Fos73}.
\end{proof}

The classes $[E_i]$ need not be linearly independent in $\Cl(R)$. The same closure argument applies to any nonempty open subset $U$ of a normal integral variety $X$: $\Cl(U)$ is the quotient of $\Cl(X)$ by the subgroup generated by prime divisors in $X\setminus U$. In particular it applies to the affine big cell of a Schubert variety even though its complement need not be the divisor of a single global function.

\begin{lemma}[Flat base change and reflection of invertibility]\label{alg:flat}
Let $A\to B$ be a flat homomorphism of Noetherian rings and let $M$ be a finite $A$-module. There are natural isomorphisms
\[
 \begin{aligned}
 \Hom_A(M,A)\otimes_A B&\simeq\Hom_B(M\otimes_A B,B),\\
 M^{**}\otimes_A B&\simeq(M\otimes_A B)^{**}.
 \end{aligned}
\]
For normal domains, flat base change therefore preserves finite rank-one reflexive modules and their reflexive tensor powers. Invertibility descends under faithful flatness; over a local source, freeness of rank one descends.
\end{lemma}
\begin{proof}
A finite presentation expresses the dual as a kernel between finite free modules. Flat tensor product preserves this kernel \cite[Tag 0C6I]{The26}. The dual is again finite over a Noetherian ring, so applying the argument twice proves the double-dual formula, and applying it to tensor powers proves the corresponding assertion. A rank-one module remains rank one after the indicated flat extension of domains. The last assertion is faithfully flat descent for finite locally free modules \cite[Tag 05B2]{The26}.
\end{proof}

\begin{lemma}[Generic smooth factors]\label{alg:generic}
Let $(R,\mathfrak m)$ be a normal Noetherian local domain and set $R(T)=R[T]_{\mathfrak mR[T]}$. A finite rank-one reflexive $R$-module $M$ is free if and only if $M\otimes_RR(T)$ is free. Consequently passage to the generic point of an affine cell factor reflects the Cartier condition. The analogous multivariable localization is denoted by $R(T_1,\ldots,T_a)$, and the same assertion holds for it.
\end{lemma}
\begin{proof}
This is a faithfully flat local extension with residue field extension $\kappa(\mathfrak m)\subset\kappa(\mathfrak m)(T)$. The minimal number of generators is unchanged. If the extended module is free of rank one, $M$ is cyclic, and a cyclic torsion-free module of rank one over a domain is free. The forward implication is immediate.
\end{proof}

\subsection{Charts and conventions}\label{sec:charts}
Write $W$ for the Weyl group, $C_u=BuB/B$ for an ordinary cell, and $C^u=B^-uB/B$ for an opposite cell. For $u\le w$, the Kazhdan--Lusztig slice $N_{u,w}=X_w\cap C^u$ is affine, has distinguished point $u$, and has dimension $\ell(w)-\ell(u)$. A standard open neighbourhood $\Omega_u$ gives an actual scheme isomorphism
\[
 X_w\cap\Omega_u\simeq C_u\times N_{u,w},
\]
compatible with Schubert strata \cite[Lemma~3.2]{WY08}\cite[Lemma~2.1]{KWY13}. Normality and Cohen--Macaulayness of the slices follow from this product and the corresponding results for Schubert varieties \cite[Theorems~2.1.1, 2.2.3 and Proposition~2.2.5]{Bri05}. The opposite cell supplies a contracting positive grading with degree-zero part $\C$.

When identifying prime ideals, we use scheme-theoretic reducedness of intersections of Schubert varieties and their opposite analogues \cite{Ram87}\cite[Section~2]{Kum17}. A calculation of root tangent directions by itself would not establish the required ideal equality.

For a global Weil divisor $D$, its non-Cartier locus $\NC(D)$ is the closed nonfree locus of the coherent reflexive sheaf $\OO_X(D)$. If $D$ is $B$-stable, so is this locus; its irreducible components are then Schubert subvarieties, because there are finitely many $B$-orbits. Any local Weil class extends to a global Weil class by closing its prime divisors. It can also be represented by a $B$-stable boundary divisor: the big cell is affine space, its class group is zero, and Lemma~\ref{alg:localization} supplies the representation.

\section{Divisor detection in type A}\label{sec:type_a}

We work over $\C$ and use the convention $X_w=\overline{BwB/B}\subset
\GL_n/B$, where $B$ is upper triangular and $\dim X_w=\ell(w)$. For a
normal variety $X$ and a Weil divisor $D$, let $\NonCart_X(D)$ be its
non-Cartier locus. The following statement concerns every Weil divisor,
including divisors whose local classes have infinite order.

\begin{theorem}[Type A divisor detection]\label{ta:main}
For every $w\in S_n$ and every Weil divisor $D$ on $X_w$,
\begin{equation}\label{ta:eq:locus}
 \NonCart_{X_w}(D)=
 \bigcup_{\substack{S\in\Irr\Sing(X_w)\\
                   D\text{ is not Cartier at }\eta_S}}S.
\end{equation}
Equivalently, for every scheme point $p\in X_w$, restriction induces an
injection
\begin{equation}\label{ta:eq:local_injection}
 \Cl(\OO_{X_w,p})\longrightarrow
 \bigoplus_{\substack{S\in\Irr\Sing(X_w)\\p\in S}}
 \Cl(\OO_{X_w,\eta_S}).
\end{equation}
The same statement holds for Schubert varieties in type A partial flag
varieties.
\end{theorem}

The proof constructs detectors for integer differences of divisor
coefficients. Its geometric ingredients are a filtered determinantal
slice construction and a procedure that replaces a horizontal component
of a first-line hyperplane section by a comparison between vertical
components.

\subsection{Coordinates and a first-difference localization}

Put $C_y=ByB/B$, $C^x=B^-xB/B$, $N_{x,y}=X_y\cap C^x$, and
$A_{x,y}=\C[N_{x,y}]$. For $x\leq y$, this slice is integral, normal and
Cohen--Macaulay, with dimension $\ell(y)-\ell(x)$.
The affine chart $\Omega_x=xC^{\id}$ satisfies
\[
 \Omega_x\cap X_y\simeq C_x\times N_{x,y}
 \simeq\A^{\ell(x)}\times N_{x,y},
\]
compatibly with the Schubert strata; see
\cite[Section 3 and Lemma 3.2]{WY08}. The standard normality and rational-singularity results,
and the Cohen--Macaulay property, pass to the slice through this product
\cite[Theorems 2.1.1, 2.2.3 and Proposition 2.2.5]{Bri05}.

For $y\in S_n$, set $r_y(r,k)=\#\{h\leq k:y(h)\geq r\}$.
Then $x\leq y$ exactly when $r_x(r,k)\leq r_y(r,k)$ for every $r,k$.
If $p<q$ and $x(p)<x(q)$, exchanging positions $p,q$ increases the rank
matrix by one precisely on the rectangle
$p\leq k<q,\; x(p)<r\leq x(q)$.
This is the rank-matrix test used below.

In the standard matrix for $C^x$, column $p$ has pivot $1$ in row $x(p)$.
Its free entries are $Z_{x(q),p}$ for $p<q$ and $x(p)<x(q)$, of weights
$e_{x(q)}-e_{x(p)}$. Distinct free entries have distinct root weights.
The equations of $N_{x,y}$ bound the southwest submatrix ranks by
$r_y(r,k)$.

\begin{lemma}\label{ta:tangent}
A free entry $Z_{x(q),p}$ has nonzero image in the vertex cotangent space
if and only if $x(pq)\leq y$. The entries with this property generate the
actual coordinate ring $A_{x,y}$.
\end{lemma}

\begin{proof}
On the corresponding coordinate axis, each defining minor is at most
linear in the free entry. The axis satisfies the rank bounds exactly when
$x(pq)\leq y$, by the rectangle test. Otherwise one defining minor
restricts to a nonzero linear term. Since the entries have distinct root
weights, no linear combination of different entries changes this criterion.

Choose a cocharacter positive on all the roots $e_{x(q)}-e_{x(p)}$ in
question. The coordinate ring is positively graded with degree-zero part
$\C$. An entry with zero cotangent class is a sum of products of elements
of smaller positive degree. Induction on degree therefore expresses every
entry as a polynomial in the surviving entries. Taking the appropriate
torus-weight component makes the expression weight-homogeneous. This
argument takes place in the coordinate ring itself.
\end{proof}

In particular, an entry whose class in the vertex cotangent space
vanishes is zero in the coordinate ring if its root weight is not a
sum of at least two surviving root weights.
If $x(1)=y(1)$, the first-column rank bounds force that column to equal
its pivot column, and its pivot row has no other nonzero entries.
Deleting that row and column and standardizing gives an actual smaller
slice, compatibly with its Schubert strata. We may thus delete a common
prefix and assume, whenever $v<w$, that $a=v(1)<w(1)=c$.
Define
\[
 j=\min\{h>1:a<v(h)\leq c\},\qquad
 b=v(j),\quad u=v(1j),\quad t=Z_{b,1}.
\]

\begin{lemma}\label{ta:open}
The permutation $u$ satisfies $v\lessdot u\leq w$, and there is an actual
isomorphism
\[
 N_{v,w}[t^{-1}]\simeq\Gm\times N_{u,w}
\]
that preserves ordinary Schubert strata.
\end{lemma}

\begin{proof}
The choice of $j$ excludes values strictly between $a,b$ in the intervening
positions, so the exchange is a Bruhat cover. On its rectangle
$1\leq k<j$, $a<r\leq b$, the rank of $v$ counts only values greater than
$c$. The bound for $w$ contains at least those high values, by the
inequality at threshold $c+1$, as well as its first value $c$. This gives
one unit of slack and proves $u\leq w$.

More generally, for a cover $u=v(ij)$, the prefix Pl\"ucker coordinates
that define $\Omega_u$ restrict on $C^v$ to $\pm t$ for $i\leq k<j$ and
to $\pm1$ otherwise, where $t=Z_{v(j),i}$. Any other nonconstant term in
the determinant expansion would require an increasing path through an
intermediate position with value strictly between $v(i),v(j)$. The empty
rectangle of a cover excludes such a path. Consequently
$N_{v,w}\cap\Omega_u=D(t)$.

The two-sided chart isomorphism of \cite[Lemma 2.1]{KWY13} is
$\Omega_u\simeq C^u\times C_u$. Its first projection preserves the
ordinary cell of a point; its second preserves the opposite cell.
The conditions of belonging to $X_w$ and to $C^v$ therefore become,
simultaneously, conditions on the first and second factors. Restriction
of that isomorphism gives
$D(t)\simeq N_{u,w}\times(C_u\cap C^v)$.
Both sides carry their reduced scheme structures. For a cover, the second
factor is the open Richardson curve $\Gm$. This proves
the asserted isomorphism without cancelling affine factors.
\end{proof}

\subsection{The reduced first-line fibre}

Let $p_r$ be the $r$th coordinate section of the pullback of
$\OO_{\PP^{n-1}}(1)$ under the map $F_\bullet\mapsto F_1$.
In the normalized chart at $v$, $p_a=1$ and $t=p_b/p_a$.
In this cell-product chart, $p_b$ is independent of the cell factor.
Indeed, the only possible additional terms in $p_b(gM)$ have the form
$g_{bq}Z_{q,1}$, where $q>b$ and $v^{-1}(q)<j$. If $q\leq c$ this
contradicts the definition of $j$; if $q>c$, the first-line rank bound
already gives $Z_{q,1}=0$. Thus $p_b(gM)=Z_{b,1}$.

Scan the values of $w$ after its first position that are less than $c$.
Write their successive record highs as $d_1<\cdots<d_m<c$, with
$w(j_s)=d_s$ and $j_1<\cdots<j_m$. The exchanges changing the first
position in a coatom are exactly
$z_s=w(1j_s)\lessdot w$.
Their rectangles are $1\leq k<j_s$, $d_s<r\leq c$.
For $b<c$, put $H_b=\overline{C_w\cap\{p_b=0\}}\subset X_w$.

\begin{lemma}\label{ta:fibre}
The parameter divisor on $N_{v,w}$ is reduced. Its vertical components
are exactly $N_{v,z_s}$ with $v\leq z_s$ and $d_s<b$. If $b<c$, its
only other possible component is the restriction of $H_b$, which occurs
exactly when $H_b$ meets $\Omega_v$. Globally,
\[
 \begin{aligned}
  \operatorname{div}(p_b)&=H_b+\sum_{d_s<b}X_{z_s}
     &&\text{for }b<c,\\
  \operatorname{div}(p_c)&=\sum_sX_{z_s}.
 \end{aligned}
\]
The resulting divisor-class relation is
$[H_b]=\sum_{d_s\geq b}[X_{z_s}]$.
In particular, $\operatorname{div}_{N_{v,w}}(t)=\sum_{E\subset V(t)}E$,
where the sum is over its prime components, each with coefficient one.
\end{lemma}

\begin{proof}
On $C_y$, the section $p_b$ is identically zero when $y(1)<b$, nowhere
zero when $y(1)=b$, and an affine coordinate when $y(1)>b$.
The first case contributes a divisor only for a coatom, and the last
only for the big cell $C_w$. The latter zero hyperplane is irreducible.
This proves the stated component lists for the two global zero divisors.

The big-cell component has multiplicity one. At the generic point of a
vertical coatom, choose a transverse root-line parameter $\tau$ so that
the first column has the form $g(e_{d_s}+\tau e_c)$. For $d_s<b\leq c$, its $b$th
entry is $g_{bc}\tau$. If $b=c$, $g_{cc}$ is a unit. If $b<c$, the
record condition permits the independent generic cell coordinate
$g_{bc}$, so this coefficient is nonzero at the generic point. Every
vertical multiplicity is therefore one.

A nonzero section on a Cohen--Macaulay variety defines a divisor with
no embedded associated primes. Generic reducedness consequently implies
reducedness. Cell-factor independence descends the assertion to the
slice. Finally $p_b,p_c$ are sections of the same line bundle; subtracting
their divisors gives the stated relation for $[H_b]$.
\end{proof}

\subsection{Determinantal divisor classes}\label{ta:determinantal}

We use Nagata localization and polynomial invariance from
Lemma~\ref{alg:localization}, and the faithful generic-cell comparison
of Lemma~\ref{alg:generic}. The following standard calculation makes
explicit the integral local assertion needed for the detectors
\cite{Fos73,BV88}.

For positive integers $r,s$, let $D_{r,s}$ be the variety of $r\times s$
matrices of rank at most one. For $r,s\geq2$, its coordinate ring
$R=\C[x_{ij}]/I_2(x)$ is the normal Segre semigroup ring of dimension
$r+s-1$. The ruling primes $P=(x_{1j})_j$ and $Q=(x_{i1})_i$ satisfy
$\operatorname{div}(x_{11})=P+Q$. Localizing gives
\[
 R[x_{11}^{-1}]
 =\C[x_{11}^{\pm1},x_{i1}\ (i>1),x_{1j}\ (j>1)],
\]
whose units are $\C^*x_{11}^{\Z}$. The divisor localization sequence,
including its map from units, gives
\[
 \Cl(R)=(\Z[P]\oplus\Z[Q])/\Z(P+Q)\simeq\Z.
\]
For completeness, normality follows because the semigroup is the set
of pairs of nonnegative integral vectors having equal coordinate sums,
which is saturated in its group.

The same class group occurs at the vertex $\mathfrak m=(x_{ij})$.
If $n[P]$ vanishes locally for $n\neq0$, set $a=|n|$.
The homogeneous divisorial ideal $P^{(a)}$ becomes principal at
$\mathfrak m$, so $P^{(a)}/\mathfrak mP^{(a)}$ is one-dimensional.
A homogeneous lift and graded Nakayama imply that it is principal
already over $R$, contradicting the global class-group calculation.
Since localization on class groups is surjective, we obtain
\[
 \Cl(R_{\mathfrak m})\simeq\Z,\qquad [P]=1,\quad [Q]=-1.
\]
Every nonzero matrix has a nonzero-entry chart that is smooth, whereas
the vertex has tangent dimension $rs>r+s-1$. Thus
$\Sing(D_{r,s})=\{0\}$ for $r,s\geq2$.
The calculation remains valid over every field $F$ in place of $\C$:
the localized ring has units $F^*x_{11}^{\Z}$ and the graded Nakayama
argument is unchanged. In particular, this local calculation applies
to the ordinary local ring over the function field of a Schubert cell.

\subsection{Filtered maxima and actual determinantal slices}

Fix records $s<t$ and a threshold $d_t<b\leq c$.
For $y$ with $y(1)<c$, let $B(y)$ denote the first value in
$(y(1),c]$ encountered after position one. Let $\mathcal F$ be the finite
set of permutations $x\leq z_s,z_t$ such that $x\nleq z_h$ for every
$s<h<t$ and $B(x)\geq b$.

\begin{theorem}[Filtered determinantal slices]\label{ta:filtered}
Suppose $\mathcal F$ is nonempty, and let $x$ be any Bruhat-maximal
element of it. There are $p,q\geq1$ and $h\geq0$ such that
$N_{x,z_t}\simeq D_{q,p}\times\A^h$, and, as actual affine varieties,
\begin{equation}\label{ta:eq:ambient_models}
 N_{x,w}\simeq
 \begin{cases}
  D_{q+1,p},&h=0,\quad p\geq2,\\
  D_{2,p+h},&h>0,\quad q=1.
 \end{cases}
\end{equation}
In particular, $X_x$ is an irreducible component of $\Sing(X_w)$.
The first-difference parameter at $x$ becomes a matrix entry, and
$X_{z_s},X_{z_t}$ become its two primitive ruling divisors.
\end{theorem}

\begin{proof}
We supply the coordinate and rank calculations because the placement
of the two specified divisors is needed, in addition to the singularity
type.

\smallskip\noindent\emph{Normalization of the maximum.}
Write $a=d_s$ and $d=d_t$. Obtain $x_0$ from $w$ by replacing the values
$(c,a,d)$ in positions $(1,j_s,j_t)$ with $(a,d,c)$. Counting separately
before $j_s$, between $j_s,j_t$, and after $j_t$, gives
$r_{x_0}=\min(r_{z_s},r_{z_t})$. Consequently,
$y\leq z_s,z_t$ if and only if $y\leq x_0$.
For $k<j_s$, the minimum is $r_w(r,k)-\mathbf1_{a<r\leq c}$;
for $j_s\leq k<j_t$, it is $r_w(r,k)-\mathbf1_{d<r\leq c}$;
and for $k\geq j_t$ it is $r_w(r,k)$. This verifies the equality
at every threshold, including those outside the rectangles.

If $x(1)<a$, the first-difference cover for $[x,x_0]$ raises the
first value to the earliest encountered value in $(x(1),a]$.
Since $B(x)\geq b>d>a$, that value occurs after $B(x)$, so the exchange
leaves $B$ unchanged. Exclusion from intermediate coatoms is preserved
on moving upward. The resulting larger element would belong to
$\mathcal F$, contrary to maximality. Hence $x(1)=a$.

Suppose $B=B(x)<c$, at position $j$, and let $k>j$ be the first
position with $Q=x(k)\in(B,c]$. The exchange of $j,k$ is an upward
cover with rectangle $j\leq u<k$, $B<r\leq Q$. On this rectangle,
the rank of $x$ counts only values greater than $c$. Denote their
prefix counts in $x,x_0$ by $V_u,W_u$. For $u<j_t$, recordhood says
that all values of $x_0$ at most $c$ in this prefix are at most
$d<B$. At threshold $B$, the prefix of $x$ already contains $B$;
thus $x\leq x_0$ gives $V_u+1\leq W_u$. For $u\geq j_t$, the
prefix of $x_0$ contains $c$, and the comparison at threshold $c+1$
gives $V_u\leq W_u$. In either case there is one unit of slack
throughout the exchange rectangle. The exchange stays below $x_0$,
raises $B$ to $Q$, and preserves all intermediate-coatom exclusions,
again contradicting maximality. We have shown $x(1)=a$ and $B(x)=c$.
If $j=x^{-1}(c)$, this means
$x(k)\leq a\ \text{or}\ x(k)>c\;(1<k<j)$.

We also need the following consequence of recordhood. Every $y\leq w$
with $y(1)=a$ lies below $z_s$. On the coatom rectangle for $z_s$,
all values of $w$ before $j_s$ other than $c$ are less than $a$ or
greater than $c$, so $r_w(a,k)=r_w(r,k)$. Since $y(1)=a<r$,
we have
\[
 r_y(r,k)+1\leq r_y(a,k)\leq r_w(a,k)=r_w(r,k).
\]
Outside the rectangle use $y\leq w$. Maximality of $x$ consequently
excludes the simultaneous conditions
$y>x,\; y\leq z_t,\; y(1)=a,\; B(y)=c$.

\smallskip\noindent\emph{The actual support of the lower slice.}
Put $z=z_t$. A surviving tangent coordinate of $N_{x,z}$ exchanges
positions $u<v$ with $x(u)<x(v)$. If $u=1$, the new first value
cannot exceed $z(1)=d<c$; by the prefix condition its target occurs
after $j$ and lies in $(a,d]$. If $u>1$, maximality
forces the exchange to destroy $B=c$. Exchanges wholly before or after
$j$ do not do so, and neither does an exchange with $v=j$.
An exchange with $u<j<v$ must bring a value of $(a,c)$ before $j$;
hence $x(u)<a<x(v)<c$. An exchange with $u=j$ must have $x(v)>c$,
and an intermediate value in $(a,c)$ must occur between $j,v$;
otherwise the new first eligible value would still be $c$.

All surviving roots therefore have the forms $L\longrightarrow M$
and $c\longrightarrow H$, where
$L\subseteq[1,a]$, $M\subseteq[a+1,c-1]$, and
$H\subseteq[c+1,n]$. The values in $L$ occur before $j$, while those
in $M,H$ occur after $j$. The first group is nonempty and $a\in L$,
by the first-difference cover for $[x,z]$.
The negative ends of these roots belong to $L\cup\{c\}$ and their
positive ends to the disjoint set $M\cup H$. No sum of at least two
of them is a root. Lemma~\ref{ta:tangent} thus makes every other raw
entry zero in the actual ring. The matrix has only its permutation
pivots, a possibly sparse block $Z$ with rows $M$ and source columns
$L$, and entries $Y_h$ in column $j$, with rows in $H$.

\smallskip\noindent\emph{The core slack.}
For prefixes $k<j$, row $c$ is zero, so ranks at thresholds $c,c+1$
agree. Evaluating on the dense ordinary stratum of the slice gives
$j\leq j_t=z^{-1}(c)$ and
$r_z(c+1,k)=r_x(c+1,k)$ for $k<j$; high rows in these prefixes have
only fixed pivots.
If an intermediate record had $j_h\leq j$, then on its coatom
rectangle the prefix condition would give
\[
 r_x(r,k)=r_x(c+1,k)\leq r_w(c+1,k)\leq r_w(r,k)-1.
\]
The final unit is contributed by $w(1)=c$. This would imply
$x\leq z_h$, contrary to the filter. Thus $j<j_h$ for $s<h<t$.
Together with $j\leq j_t$ and the record property of $a$, this says
that $w(2),\ldots,w(j-1)$ contain no value in $(a,c]$.
In the same prefix of $z$, the sole value in $(a,c]$ is its first
value $d$. Consequently
$r_z(a+1,j-1)-r_x(a+1,j-1)=1$.
The southwest block with rows at least $a+1$ and columns at most
$j-1$ contains all of $Z$. Its source pivot rows lie in $L$ and are
absent; its target pivot columns lie after $j$ and are absent.
Every other nonzero fixed pivot has its own row and column disjoint
from $Z$. Removing these pivots from the rank condition leaves
exactly $\rank Z\leq1$.

\smallskip\noindent\emph{Recovering the whole lower slice.}
Set $p=|L|$, $q=|M|$, and $h=|H|$. Every active row and column of
$Z$ has a nonzero entry. If $Z_{ml}$ were missing, choose nonzero
entries $Z_{m'l}$ and $Z_{ml'}$. The corresponding minor would give
$Z_{m'l}Z_{ml'}=0$, contrary to integrality. Thus $Z$ is a complete
$q\times p$ matrix.

We have a closed torus-invariant embedding
$N_{x,z}\subseteq D_{q,p}\times\A^h$. The product of its surviving
coordinates is nonzero in its domain, so the slice has a point where
every $Z,Y$ entry is nonzero. All such points of the ambient product
form one torus orbit: write $Z_{ml}=u_mv_l$ and $Y_h=y_h$, with all
factors nonzero, and choose diagonal entries
\[
 t_l=v_l^{-1},\qquad t_m=u_m,\qquad t_c=1,\qquad t_h=y_h.
\]
The disjoint value sets make the choices compatible. This orbit is
dense, so a closed invariant subvariety containing one of its points
contains its closure. Therefore
$N_{x,z}\simeq D_{q,p}\times\A^h$; since $z\lessdot w$,
$\dim N_{x,w}=p+q+h$.

\smallskip\noindent\emph{A dense torus orbit in the ambient slice.}
The surviving ambient roots include the lower roots and the
first-difference root $a\longrightarrow c$. Their graph connects
$V=L\cup M\cup\{c\}\cup H$ and has root rank
$|V|-1=p+q+h$. This is the usual rank of the incidence vectors
$e_v-e_u$ of a connected graph. Choose a point $P\in N_{x,w}$ where
all surviving coordinates are nonzero. Its torus orbit has dimension
equal to the rank of the surviving weights. This dimension is at least
$p+q+h$ and at most $\dim N_{x,w}=p+q+h$; hence the orbit is dense.

Normalize each surviving coordinate $f_\alpha$ by its nonzero value
at $P$. On the orbit it then restricts to the character $\chi^\alpha$.
Because these entries generate the actual coordinate ring, we have
\[
 \C[N_{x,w}]=\C[\mathbb N\{\alpha:\alpha\text{ survives}\}].
\]
Indeed, restriction to the dense orbit is injective and its image is
exactly the algebra generated by these normalized characters.
Equal character monomials are therefore equal in the actual slice
ring.

An additional surviving root cannot involve a value outside
$V$, since that would increase root rank. Nor can a surviving root
be a sum of two others, since its normalized coordinate would then
belong to the square of the vertex ideal. The normalization multiplies
coordinates only by nonzero constants and therefore preserves their
zero ideals, including the marked divisors used below.

In the following southwest-rank calculation we retain the original
lower-slice matrix entries with their fixed pivots. The normalized
characters are used when writing the final ambient determinantal equations.

\smallskip\noindent\emph{Generic southwest pivots.}
Index the source positions by $\ell_1=1<\cdots<\ell_p<j$, the
middle positions by $j<m_1<\cdots<m_q$, and the high positions by
$j<h_1<\cdots<h_h$. Write $x(\ell_i)=\lambda_i$, $x(m_r)=\mu_r$,
and $x(h_i)=\eta_i$. Within each list the values strictly decrease:
\[
 \lambda_1=a>\cdots>\lambda_p,\qquad
 \mu_1>\cdots>\mu_q,\qquad
 \eta_1>\cdots>\eta_h>c.
\]
For instance, an increasing source pair has a reflection rectangle
contained in that from the earlier source to an active middle target.
It would give an exchange below $z$ preserving the first value and
$B=c$, contrary to maximality. The first source cannot occur
in such an increasing pair, because its value is $a$. For middle or high targets,
use an edge to the later target to contain the increasing-pair
rectangle; the exchange again preserves the same data.

For descending $\lambda_i,\mu_r$, with every $\mu_r>\lambda_i$,
consider successive source columns
$e_{\lambda_i}+v_iU$, where $U=\sum_{r=1}^q u_re_{\mu_r}$ and
$u_rv_i\neq0$, followed by target columns
$e_{\mu_1},\ldots,e_{\mu_q}$. For a prefix of $t\geq1$ source
columns, if a threshold retains $A$ source rows, its rank is
$A+\mathbf1_{\mu_1\geq r,\ A<t}$. These are the ranks of the
values $\mu_1,\lambda_1,\ldots,\lambda_{t-1}$.

After all source columns and $u<q$ target columns, if a threshold
retains $A$ source rows and $B$ target rows, its rank is
$A+\min(u,B)+\mathbf1_{A<p,\ B>u}$. These are the ranks of
$\mu_1,\ldots,\mu_{u+1},\lambda_1,\ldots,\lambda_{p-1}$.

The last target restores all original pivot rows. The generic pivot
sequence is therefore
$\mu_1,\lambda_1,\ldots,\lambda_{p-1},
\mu_2,\ldots,\mu_q,\lambda_p$.
Applying this to the dense rank-one locus of $Z$ gives
\begin{align}
 z(\ell_1)&=\mu_1,&
 z(\ell_i)&=\lambda_{i-1}\quad(2\leq i\leq p),\notag\\
 z(m_r)&=\mu_{r+1}\quad(r<q),&
 z(m_q)&=\lambda_p.\label{ta:eq:middle_pivots}
\end{align}
The star is the one-source case of the same calculation: when $h>0$,
$z(j)=\eta_1$, $z(h_i)=\eta_{i+1}$ for $i<h$, and $z(h_h)=c$.
If $h=0$, $z(j)=c$. Other positions keep their $x$-values.
Interlacing the two groups does not change the calculation because
their row and column supports are disjoint: every southwest rank is
the sum of their ranks and those of the independent fixed pivots.
The dense ordinary stratum is indexed by $z$, so these are its actual
values. In particular, $\mu_1=d$, while $j_t=j$ for $h=0$ and
$j_t=h_h$ for $h>0$. The permutation $w=z(1j_t)$ exchanges the first
value $d$ with $c$.

\smallskip\noindent\emph{The possible parameters.}
Put $K=\min(m_1,j_t)$. We claim that
\[
 r_w(r,k)-r_x(r,k)\geq1
 \qquad\text{whenever }a<r\leq c\text{ and }1\leq k<K.
\]
For $k<j$, use the high-threshold comparison and $w(1)=c$.
For $j\leq k<K$, all sources but no middle target have appeared.
Formula \eqref{ta:eq:middle_pivots} changes the ranks in $(a,c]$
by $\mathbf1_{r\leq d}$; the star does not change these counts;
and $z\to w$ contributes $\mathbf1_{d<r\leq c}$. The sum is one.
If an intermediate record occurred at a position at most $K$, its
entire coatom rectangle would be covered by
the rank bound just proved, giving $x\leq z_h$.
Thus no intermediate record occurs at or before $K$.

If $h=0$ and $p=1$, \eqref{ta:eq:middle_pivots} puts $a$ at position
$m_q>j=j_t$ in $w$, contrary to $j_s<j_t$. Hence $p\geq2$ in this
case. Now suppose $h>0$ and $q\geq2$. If $m_1<j_t$, then
$w(m_1)=\mu_2\in(a,d)$. If $a$ has not yet occurred in $w$, it
cannot subsequently be a record. Otherwise a record strictly between
$a,d$ has appeared at or before $m_1=K$, again a contradiction.
Thus $m_1>j_t$.

The surviving star edge $c\to\eta_1$ must have a value in $(a,c)$
between positions $j,h_1$; without it the exchange would preserve
$B=c$ and contradict maximality. All active middle
targets lie after $m_1>j_t\geq h_1$, so this intermediate value
belongs to neither active group and is unchanged in $w$.
It occurs before the record $d$, so it is less than $d$.
It either occurs before $a$, destroying recordhood of $a$, or forces
an intermediate record before $j_t=K$. Both are impossible; hence
$q=1$. With a unique middle position $m_1$, the same argument excludes
$m_1>h_1$: the required intermediate value would again be unchanged
and occur before $\min(m_1,j_t)=K$. We have proved that $p\geq2$
when $h=0$, while $q=1$ and $m_1<h_1$ when $h>0$.

\smallskip\noindent\emph{All ambient roots and the actual equations.}
Every $l\in L$ has an ambient edge $l\to c$. For $l=a$ this is
the first-difference edge. For $l<a$, the known rectangles of
$l\to d$ and $a\to c$ cover that of $l\to c$, because
$l<a<d<c$ and $\operatorname{pos}(l)<j<\operatorname{pos}(d)$.
There are no roots outside $V$, by the rank argument, and no
increasing roots within $L,M,H$, by the strict decrease in each value list.
No $L\to H$ root survives: it is the sum of $l\to c$ and
$c\to h$, contradicting the semigroup-algebra identity above and tangent survival.
Only $M\to H$ remains to be considered. Thus for $h=0$ the full
ambient list is $L\to M\cup\{c\}$.

For $h>0$, we have $M=\{d\}$ and its position precedes every high
target. Every $d\to\eta_i$ survives. Above threshold $c$, its
rectangle is contained in that of $c\to\eta_i$.
For $d<r\leq c$, the prefix has passed the sole middle target, so
the lower middle block has restored its original ranks, and the star
does not change these counts. The exchange $z\to w$ adds one
because $\operatorname{pos}(d)\leq k<h_i\leq j_t$.
The complete ambient list is therefore $L\to\{d,c\}\to H$.

In the case $h=0$, arrange the normalized characters in a matrix with
rows $M\cup\{c\}$ and columns $L$. Opposite products in each
$2\times2$ minor have the same character, so all such minors vanish.
Its entries generate the ring, giving a closed embedding into
$D_{q+1,p}$. Both varieties are integral of dimension $p+q$, so
this is an isomorphism. When $h>0$, arrange them as
\[
 \begin{pmatrix}
 X_{c,l}\ (l\in L)&Y_{h,d}\ (h\in H)\\
 X_{d,l}\ (l\in L)&Y_{h,c}\ (h\in H)
 \end{pmatrix}.
\]
The difference between its two row weights is always $e_c-e_d$,
so the same character argument makes its minors zero. All generators
occur and both dimensions are $p+h+1$, giving $D_{2,p+h}$.
This proves \eqref{ta:eq:ambient_models}; no affine factor has been
suppressed.

The parameter restrictions just established make both matrix dimensions at least two.
The ambient slice is singular exactly at its vertex, so the actual
cell chart gives $\Sing(X_w)\cap\Omega_x=C_x$.
The singular locus is closed and $B$-stable and therefore contains
$X_x$. A larger singular Schubert component containing $X_x$ would
have a nonempty open intersection with $\Omega_x$ of dimension
larger than $\dim C_x$, a contradiction. Thus $X_x$ is a maximal
singular component.

Finally, the prefix parameter $Z_{c,1}$ is the surviving root
$a\to c$, hence a nonzero scalar times the indicated matrix entry.
Its zero divisor has exactly two reduced ruling primes.
The distinct coatom slices $N_{x,z_s},N_{x,z_t}$ are codimension-one
subvarieties in this zero divisor because their first values are less
than $c$, so they are precisely the two rulings.
The local determinantal calculation gives their primitive opposite classes.
\end{proof}

Two marked models illustrate the cases. For $w=4231$, take the records
$a=2$, $d=3$ and threshold $b=c=4$. The filtered maximum is $x=2143$,
with $z_s=2431$, $z_t=3241$, and $(p,q,h)=(2,1,0)$.
The opposite-cell matrix and remaining equation are
\[
 \begin{pmatrix}
 0&1&0&0\\1&0&0&0\\B&D&0&1\\A&C&1&0
 \end{pmatrix},\qquad AD-BC=0.
\]
The determinantal matrix is $\bigl(\begin{smallmatrix}A&C\\B&D\end{smallmatrix}\bigr)$,
the parameter is $A$, and the two prescribed coatom ideals are
$(A,B)$ and $(A,C)$.

For the other case, take $w=3412$, $a=1$,
$d=2$, $b=c=3$. The filtered maximum is $x=1324$, with $z_s=1432$,
$z_t=2413$, and $(p,q,h)=(1,1,1)$. Its matrix is
\[
 \begin{pmatrix}
 1&0&0&0\\B&0&1&0\\A&1&0&0\\C&D&E&1
 \end{pmatrix}.
\]
The southwest rank ideal is $(C,AD+BE)$, so the actual slice has
matrix $\bigl(\begin{smallmatrix}B&A\\-D&E\end{smallmatrix}\bigr)$ of rank
at most one. Again the parameter is $A$, and the marked ideals are
$(A,B)$ and $(A,E)$. In both examples the southwest rank conditions
are exactly the displayed equations; the rank criterion also checks
the filtered maximum directly.

\subsection{Removing a horizontal component}

Return to $a=v(1)<b=v(j)\leq c=w(1)$. The vertical components of
its fibre are the active records $v\leq z_s$ with $d_s<b$.
If there are any, let $z_*=w(1h_*)$ be the one with largest first
value $d=z_*(1)<b$. We retain its failure to be Cartier while
moving to a point whose first eligible value is $c$.

\begin{lemma}[Horizontal incidence]\label{ta:incidence}
For $v\leq w$ and $b<c$, whether or not $b=B(v)$, set
$j_b=v^{-1}(b)$ and $L_b=\min\{h>1:b\leq w(h)<c\}$.
Let $V_k=\#\{h\leq k:v(h)>c\}$ and
$W_k=\#\{h\leq k:w(h)>c\}$. Then
\[
 H_b\cap\Omega_v\neq\varnothing
 \quad\Longleftrightarrow\quad
 V_k<W_k\text{ for every }j_b\leq k<L_b.
\]
An empty range imposes no condition. In particular, exchanging two
values less than $b$, while remaining below $w$, preserves incidence.
\end{lemma}

\begin{proof}
For $g\in B$, the first line of $gwB$ is generated by $ge_c$.
The inverse image of $C_w\cap H_b$ is the irreducible hypersurface
$g_{bc}=0$ in $B$. The point $gwB$ lies in $\Omega_v$ exactly when
every prefix minor $\det g_{I_k,J_k}$ is nonzero, where $I_k=v([k])$
and $J_k=w([k])$.
The intersection $\Omega_v\cap H_b$ is nonempty exactly when it meets the
dense subset $H_b\cap C_w$. On an irreducible hypersurface, finitely
many minors are simultaneously nonzero somewhere exactly when none
vanishes identically there.

For $I=I_k$, $J=J_k$, determinant monomials correspond to perfect
matchings in the bipartite graph with row vertices $I$, column vertices
$J$, and edges $r\leq s$. The entries, including invertible diagonal
entries, are independent, so distinct matching monomials do not cancel.
There is a matching because $v\leq w$. Deleting the edge $(b,c)$
can violate a Hall inequality only for a row subset $S$ containing
$b$ that loses column $c$ from its neighborhood. Every other row of
$S$ must exceed $c$. Its neighborhood after deletion is then
$(J\cap[b,n])\setminus\{c\}$, independent of those other rows.
The largest possible $S$ is $\{b\}\cup(I\cap(c,n])$.
Since the original Hall inequalities hold, failure occurs precisely
when $b\in I$ and $1+\#(I>c)=\#(J\cap[b,n])$.
Because $c\in J$ and $\#(I>c)\leq\#(J>c)$, this is equivalent to
\[
 b\in I,\qquad J\cap[b,c)=\varnothing,\qquad
 \#(I>c)=\#(J>c).
\]
These conditions say $k\geq j_b$, $k<L_b$, and $V_k=W_k$.
Their absence is precisely the claimed incidence criterion.
A swap of two values less than $b$ changes neither $j_b$ nor the
high-value counts.
\end{proof}

\begin{lemma}\label{ta:two_components}
If the first-difference fibre has both a horizontal and a vertical
component, there is $x$ with $v\leq x\leq z_*$, $x(1)=d$, and
$B(x)=b$ whose fibre has precisely the components $z_*$ and $H_b$.
The divisor $z_*$ is not Cartier at the slice vertex.
\end{lemma}

\begin{proof}
While the first value $a_0$ is less than $d$, exchange it with the
earliest value in $(a_0,d]$. The first-difference cover for the upper
bound $z_*$ keeps the new permutation below $z_*$. Since $b>d$ is
the first value in $(a_0,c]$, the exchanged value occurs after $b$.
The exchange preserves the first eligible value $b$ and involves two
values less than $b$; Lemma~\ref{ta:incidence} preserves horizontal
incidence. Lemma~\ref{ta:fibre} applies at each step.
The first value strictly increases until it reaches $d$.
Lower records are now excluded by the first-column Bruhat bound;
records in $(d,b)$ were initially inactive and cannot become active
on moving upward. The sole remaining vertical component is $z_*$.

In a normal Noetherian local domain, if
$t\in\mathfrak m\setminus\mathfrak m^2$ and
$\operatorname{div}(t)=E_1+\cdots+E_r$ is reduced, no nonempty proper
subsum is Cartier. Otherwise that effective subsum is
$\operatorname{div}(f)$ with $f\in\mathfrak m$ by normality; its
effective nonempty complement is $\operatorname{div}(t/f)$ with
$t/f\in\mathfrak m$, forcing $t\in\mathfrak m^2$.
Here the parameter survives in the tangent space by
Lemma~\ref{ta:tangent}. Both components are invariant for the
contracting torus and contain the vertex. Applying the observation to
$\operatorname{div}(t_b)=z_*+H_b$ proves the last assertion.
\end{proof}

\begin{lemma}[The reduced first-column prime]\label{ta:column}
Let $d$ be a first-position record of $w$, let $z_*=w(1h_*)$ with
$w(h_*)=d$, and suppose $y\leq w$ and $y(1)=d$.
Then $y\leq z_*$. In $N_{y,w}$, the reduced prime ideal of
$N_{y,z_*}$ is $J_d=(Z_{p,1}:p>d)$. At the vertex its minimal
number of generators is $\mu(J_d)=q(y,w)$, where
\[
 q(y,w):=\#\{k>1:y(k)>d,\ y(1k)\leq w\}.
\]
It is Cartier there if and only if $q(y,w)=1$.
\end{lemma}

\begin{proof}
Relative to $w$, $z_*$ lowers the ranks by one on
$1\leq k<h_*$, $d<r\leq c$. If $H_k$ counts values greater than
$c$ in the corresponding prefix of $w$, recordhood gives
$r_w(d,k)=H_k+1$. Since $y(1)=d<r$,
$r_y(r,k)+1\leq r_y(d,k)\leq H_k+1$, proving $y\leq z_*$.

Consider the Schubert variety $Y_d=\{F_1\subset E_d\}$.
The support of $X_w\cap Y_d$ is the union of first-position record
coatoms with first value at most $d$. To see this, take $x\leq w$
with $a_0=x(1)\leq d$ and choose the first record $d_0\geq a_0$,
at position $h$. Such a record exists and $d_0\leq d$.
Before $h$, all nonfirst values of $w$ are less than $a_0$ or greater
than $c$. On the coatom rectangle of $w(1h)$, the inequality
$r_x(r,k)+1\leq r_x(a_0,k)\leq r_w(a_0,k)$ supplies the unit of
slack needed for $x\leq w(1h)$. The reverse support inclusion is
immediate.

We also need a scheme-theoretic assertion. Intersections of Schubert
varieties are reduced; the compatible Frobenius-splitting result and
its characteristic-zero form are recalled in
\cite[Section 2, p.~2474]{Kum17}. On $\Omega_y$, lower records are
absent because $y(1)=d$. Thus the intersection is scheme-theoretically
$X_{z_*}\cap\Omega_y$. The condition $F_1\subset E_d$ has ideal
generated by first-column entries in rows $p>d$.
Multiplication by the upper triangular cell factor acts on these tail
entries by an invertible triangular matrix and preserves their ideal.
The reduced intersection descends through the cell-product chart to the
asserted equality $J_d=(Z_{p,1}:p>d)$ of reduced prime ideals.
Its height is one.

By Lemma~\ref{ta:tangent}, the surviving first-column entries are
exactly those counted by $q(y,w)$. Their distinct weights make their
cotangent classes linearly independent. They generate all of $J_d$:
express any first-column entry of weight $e_p-e_d$ as a
weight-homogeneous polynomial in the surviving entries.
View a monomial as a nonnegative integral flow along edges from
smaller to larger values. A flow of weight $e_p-e_d$ is a path from
$d$ to $p$: inflow equals outflow elsewhere and there are no directed
cycles. It contains an edge with source $d$, which, since $y(1)=d$,
is a surviving first-column entry. Every monomial belongs to their
ideal. Their classes remain independent in $J_d/\mathfrak mJ_d$
because $\mathfrak mJ_d\subseteq\mathfrak m^2$, proving
the asserted generator count $\mu(J_d)=q(y,w)$. A height-one prime in a normal local
domain is Cartier exactly when its divisorial ideal is principal.
\end{proof}

\begin{lemma}[Greedy prefix lift]\label{ta:greedy}
Suppose $y\leq w$, $y(1)=d<c$, and $q(y,w)\geq2$.
There is $y'\geq y$ with $y'\leq w$, $y'(1)=d$, $B(y')=c$, and
$q(y',w)\geq2$. If $d$ is a record with associated coatom $z_*$,
all permutations in the construction lie below $z_*$.
\end{lemma}

\begin{proof}
If $B(y)=c$ there is nothing to prove. Otherwise put
\[
 \begin{aligned}
  j&=\min\{h>1:d<y(h)\leq c\},&
  B&=y(j)<c,\\
  k&=\min\{h>j:B<y(h)\leq c\},&
  Q&=y(k),\qquad \widetilde y=y(jk).
 \end{aligned}
\]
The position $k$ exists because $c$ occurs
after $j$. No intervening value lies strictly between $B,Q$, so this
is an upward cover with rectangle
$j\leq s<k,\; B<r\leq Q$.
On this rectangle every value of $y$ at most $c$ in the prefix is
less than $r$, so its rank counts only values above $c$.
The threshold-$c+1$ comparison with $w$, together with the first
value $c$, supplies one unit of slack. Therefore
$\widetilde y\leq w$, its first value is $d$, and its first eligible
value is $Q>B$.

Call $T>d$ a target if exchanging its position with position one
stays below $w$. Every target is at most $c$ by the first-line bound;
$B$ is a target by Lemma~\ref{ta:open}. Choose a second target
$T\neq B$. The new first-difference target is $Q$; we exhibit
another.

If $T<B$, put $h=y^{-1}(T)$ and $Y=y(1h)\leq w$.
Here $h>j$ and $h\neq k$, so $\widetilde y(1h)=Y(jk)$.
On this rectangle, all values at most $c$ in $Y$
remain below $r$, including the new first value $T$ and the value
$d$ at $h$. The high-threshold comparison for $Y\leq w$ and the
first value $c$ of $w$ again supply one unit of slack.
Thus $Y(jk)\leq w$, and $T<Q$ remains a target.

If $T>B$, put $h=y^{-1}(T)$; the definition of $k$ gives $h\geq k$.
Set $Y_T=y(1h)\leq w$, $Y_B=y(1j)\leq w$, and
$X=\widetilde y(1k)$. At positions $(1,j,k)$, the values of
$y,\widetilde y,X$ are respectively $(d,B,Q)$, $(d,Q,B)$, and
$(B,Q,d)$. For a prefix $s<j$, $X$ has the same values as $Y_B$.
For $j\leq s<k$ and $d<r\leq B$, its rank is the rank of $y$
plus one, equal to that of $Y_T$ because $h\geq k$ and $T>B$.
At all other thresholds in this middle prefix range, its rank equals
that of $\widetilde y$. For $s\geq k$, $X$ and $\widetilde y$
have the same prefix value sets. Every rank of $X$ is bounded by
that of one of $Y_B,Y_T,\widetilde y$, all below $w$.
Hence $X\leq w$, so $B<Q$ is the other target.

Each lift increases $B$, preserves $d$ and $q\geq2$, and remains
below $w$. Iteration terminates at $B=c$. When $d$ is a record,
Lemma~\ref{ta:column} puts every step below $z_*$.
\end{proof}

\begin{proposition}[Horizontal-to-prefix witness]\label{ta:horizontal}
Suppose the original first-difference fibre has both a horizontal and
a vertical component. There is $x\geq v$ with $x(1)=d$, $B(x)=c$,
and $x\leq z_*$ such that its prefix fibre has the following
properties: its sole active record of first value less than $b$ is
$z_*$; it has an active record of first value at least $b$; and
$z_*$ is not Cartier at the slice vertex.
\end{proposition}

\begin{proof}
Lemma~\ref{ta:two_components} gives a point of first value $d$ with
two-component fibre $z_*+H_b$ and with $z_*$ non-Cartier.
Lemma~\ref{ta:column} gives $q\geq2$, and Lemma~\ref{ta:greedy}
raises $B$ to $c$ while retaining this inequality and the coatom
upper bound. Thus $z_*$ remains non-Cartier at the endpoint.

All moves are upward, so inactive records cannot become active.
The first value $d$ excludes lower records, leaving $z_*$ as the
only record below $b$. The prefix fibre has no horizontal component,
by Lemma~\ref{ta:fibre}. If $z_*$ were its only component,
reducedness would give $\operatorname{div}(p_c/p_d)=z_*$, making
it Cartier, a contradiction. A further active record therefore exists,
and its first value is at least $b$.
\end{proof}

\subsection{Integral coefficient detectors and induction}

Fix $V(t)$ in $N_{v,w}$. Order its active vertical components
$E_1,\ldots,E_k$ by record first value, and place the horizontal
component $H$, if present, last. A class supported on this fibre has
a representative
\[
 D=\sum_{i=1}^k c_i[E_i]+c_H[H],\qquad c_i,c_H\in\Z,
\]
omitting the horizontal term when absent. This representation need
not be unique.

\begin{lemma}[Coefficient differences]\label{ta:differences}
For each adjacent pair $E_i,E_{i+1}$, there is a maximal singular
Schubert component $X_y$ containing $C_v$ at whose generic point the
local image of $D$ is, up to choice of a
primitive generator, $\pm(c_i-c_{i+1})\in\Z$.
If $H$ and a vertical component both occur, there is such a component
detecting $\pm(c_k-c_H)$.
\end{lemma}

\begin{proof}
Let $z_s,z_t$ correspond to adjacent verticals. At $v$, no intervening
record is active and $B(v)=b>d_t$, so $v$ belongs to
the filtered set $\mathcal F$ of Theorem~\ref{ta:filtered}. Choose a maximal element $y$ of that
set above $v$ and apply Theorem~\ref{ta:filtered}.
Its prefix fibre is the zero divisor of one matrix entry, with
exactly the rulings $z_s,z_t$. Lemma~\ref{ta:fibre} then excludes
every other active record in that chart. If the original $H_b$
exists, its class is the sum of high records $d_h\geq b$ by
the class relation in Lemma~\ref{ta:fibre}; all are absent from this chart.
The ordinary detector ring is
\[
 \OO_{X_w,\eta_{X_y}}
 \simeq(\C(C_y)[z_{ij}]/I_2(z))_{(z_{ij})},
\]
for the matrix size supplied by Theorem~\ref{ta:filtered}.
The two ruling classes are primitive opposites by the calculation in
Section~\ref{ta:determinantal}, which applies over the field $\C(C_y)$.
This gives the stated integer difference. If $P$ is the corresponding
slice prime in $A_{v,w}$ and $R=(A_{v,w})_P$, the cell-product chart also
identifies this detector with $R(T_1,\ldots,T_{\ell(v)})$.
The comparison $R\to R(T_1,\ldots,T_{\ell(v)})$ is faithfully flat and local,
so Lemma~\ref{alg:generic} reflects vanishing.

For the last difference, Proposition~\ref{ta:horizontal} gives
$x\geq v$ with first value $d$ and $B(x)=c$, where only the highest
original low record $z_*=E_k$ remains and at least one high record
is active. Choose the active high record $z_t$ of smallest first
value. It is adjacent to $z_*$ in this prefix fibre.
Use the filtered construction with threshold $d_t+1\leq c$ and a
maximal element $y\geq x$. Its two rulings are exactly $z_*,z_t$.
Every other original low record is absent, while
the class relation in Lemma~\ref{ta:fibre} restricts to
$[H_b]=[X_{z_t}]=-[X_{z_*}]$. The local coefficient of the original
class is consequently $\pm(c_k-c_H)$.
This divisor-class relation supplies the comparison without requiring
an equation for the old horizontal divisor in the new chart.
\end{proof}

\begin{theorem}[Affine slice detection]\label{ta:affine}
For every $v\leq w\in S_n$, restriction induces an injection
\begin{equation}\label{ta:eq:affine_injection}
 \Cl(N_{v,w})\longrightarrow
 \bigoplus_{S\in\Irr\Sing(N_{v,w})}
 \Cl(\OO_{N_{v,w},\eta_S}).
\end{equation}
\end{theorem}

\begin{proof}
Induct on $\ell(w)-\ell(v)$; the zero-dimensional case is immediate.
Delete a common prefix and choose the first-difference cover and
parameter from Lemma~\ref{ta:open}.
Let $D$ map to zero on the right of \eqref{ta:eq:affine_injection}.
A nonempty open subset retains exactly the irreducible singular
components meeting it; it does not turn a deeper singular stratum
into a new component. By Lemma~\ref{ta:open} and Laurent
polynomial invariance, the class on $D(t)$ comes from
$D'\in\Cl(N_{u,w})$. The maximal singular components of this
product are $\Gm\times S'$. At their generic points,
Lemma~\ref{alg:generic} shows that $D'$ vanishes at all maximal
singular generic points of $N_{u,w}$. Since its dimension is one less,
induction gives $D'=0$.

Lemma~\ref{alg:localization} writes $D$ as an integer combination of
the vertical classes $[E_i]$ and, when present, the horizontal class
$[H]$. Every detector in
Lemma~\ref{ta:differences} is a singular component of the original
slice: its global maximal singular Schubert component meets
$\Omega_v$, and its intersection is the product of $C_v$ with
its slice factor. The vanishing hypothesis therefore gives
$c_i=c_{i+1}$ for $i<k$, and $c_k=c_H$ when both kinds occur.
All coefficients equal some $m\in\Z$.
This is automatic for a single fibre component, including a sole
horizontal component. Reducedness from Lemma~\ref{ta:fibre} now gives
\[
 D=m\sum_{E\subset V(t)}[E]=m[\operatorname{div}(t)]=0.\qedhere
\]
\end{proof}

\begin{proof}[Proof of Theorem~\ref{ta:main}]
Take a scheme point $p\in C_v\subset X_w$ and
$\alpha\in\Cl(\OO_{X_w,p})$ vanishing at all maximal singular
generic points whose closures contain $p$. Put
\[
 A_{\mathrm{ch}}=\C[\Omega_v\cap X_w]
 \simeq A_{v,w}[t_1,\ldots,t_{\ell(v)}].
\]
Localization on class groups is surjective by taking closures of the
height-one primes in a representative. Choose a lift
$\beta\in\Cl(A_{\mathrm{ch}})$ of $\alpha$.

Every maximal singular component is a Schubert variety $X_z$.
It meets $\Omega_v$ exactly when $v\leq z$; in that case
the cell-product chart shows that it contains all of $C_v$, hence $p$
even if $p$ is not closed. Its generic point is a generalization
of $p$, so the restriction of $\beta$ factors through $\alpha$.
Thus every chosen lift vanishes at all maximal singular generic points
of this chart.

Polynomial invariance writes $\beta$ as the pullback of
$\gamma\in\Cl(A_{v,w})$. The full-cell generic local ring is
$R(T)$ for the corresponding slice generic local ring $R$.
Lemma~\ref{alg:generic} therefore makes $\gamma$ vanish at every
maximal singular generic point of the slice.
Theorem~\ref{ta:affine} gives $\gamma=0$, hence $\beta=0$ and
$\alpha=0$. This proves \eqref{ta:eq:local_injection} without any
equivariance assumption on $\alpha$.

The invertible locus of a finite rank-one reflexive sheaf is open.
Non-Cartierness at $\eta_S$ therefore persists on $S$, proving one
inclusion in \eqref{ta:eq:locus}. If $p$ avoids the indicated union,
the local class of $D$ vanishes at every relevant generic point;
the injection makes it zero at $p$, proving the other inclusion.

Finally, for $P\supset B$, the projection $\pi:G/B\to G/P$ is
locally trivial with irreducible smooth fibre $P/B$.
The inverse image of a partial-flag Schubert variety is the full-flag
Schubert variety indexed by the longest representative of its coset.
Its singular locus is the inverse image of the singular locus
downstairs, and inverse images of irreducible components are
irreducible. Weil divisors pull back on these local products;
Cartierness is equivalent before and after the faithfully flat smooth
extension. The generic-point maps are compatible with $\pi$.
Apply \eqref{ta:eq:locus} upstairs and descend to obtain the
partial-flag assertion.
\end{proof}

\section{Local geometry and torsion detection}\label{tor:local}

The arguments in this section apply to every finite crystallographic root system. We work over $\C$ and allow the field $F\supset\C$ to vary. This extra quantifier is necessary when a slice is replaced by a closed subslice. In any one argument, the strict henselizations of all subslices are formed using the same separable closure $F^{\mathrm{sep}}$ of $F$.

A singular label $u$ of $X_w$ means $X_u\subseteq\Sing(X_w)$; it is nonmaximal when $X_u$ is not an irreducible component of this locus. This terminology always concerns singularity of the ambient variety $X_w$ along $X_u$.

\subsection{Slices and shared local rings}\label{tor:charts}
Write $C_z=BzB/B$, $C^u=B^-uB/B$, and $N_{u,z}=X_z\cap C^u$ for $u\le z$. Intersections have their actual scheme structures. The slice has distinguished vertex $uB/B$ and dimension $\ell(z)-\ell(u)$. We use the multiplication chart
\[
 O_u\cap X_z\simeq C_u\times N_{u,z},
 \qquad O_u=u(B^-B/B).
\]
The same unipotent multiplication constructs these charts for every $z$, so they respect inclusions; see \cite[Lemma~2.1 and its proof]{KWY13}. Schubert varieties are normal and Cohen--Macaulay, and therefore their slices have these properties as well \cite[Sections~2.1--2.2]{Bri05}.

For $F\supset\C$, set
\[
 \begin{gathered}
 R_F(u,z)=F[N_{u,z}]_{\mathfrak m_u},
 \qquad A_F(u,z)=R_F(u,z)^{\sh},\\
 U_F(u,z)=\Spec A_F(u,z)\setminus\{s\},
 \end{gathered}
\]
where $\mathfrak m_u$ is the vertex ideal and $s$ is the closed point. The rings $R_F(u,z)$ and $A_F(u,z)$ are Noetherian excellent normal local domains, and they are Cohen--Macaulay. Indeed, over the algebraically closed field $\C$, normal integral varieties are geometrically normal and geometrically integral. Normality persists under extension of the ground field \cite[Tags~038N and 038P]{The26}; the Cohen--Macaulay property does so here because the algebras are of finite type \cite[Tag~045P]{The26}. Strict henselization preserves the required local properties \cite[Tags~06DI and 06LM]{The26}.

Fix $u\le z\le w$, and let $I_z$ be the ideal of the closed subslice. Compatibility of strict henselization with quotients identifies
\[
 A_F(u,z)=A_F(u,w)/I_zA_F(u,w)
\]
with the chosen common separable closure \cite[Tag~05WS]{The26}. Consequently the $U_F(u,z)$ are closed subschemes of a common punctured spectrum. Each is nonempty and connected when $z>u$.

The product chart also gives
\[
 \OO_{X_w,\eta_{X_u}}
 \simeq\bigl(\C(C_u)[N_{u,w}]\bigr)_{\mathfrak m_u}.
\]
Thus the ordinary generic-point ring is obtained by taking $F=\C(C_u)$; the smaller subslices in the same argument must still use this field.

\subsection{Reduction to an effective Levi subgroup}
\begin{proposition}\label{tor:levi}
If $N=N_{u,w}$ has dimension $d$, it is isomorphic, as an affine scheme with vertex and ambient torus action, to a Schubert slice for a Levi subgroup whose derived rank is at most $d$. The isomorphism preserves root lengths after conjugating that Levi to a standard Levi subgroup, and it commutes with every extension $F\supset\C$.
\end{proposition}
\begin{proof}
Restrict the root-coordinate functions of the opposite cell to $N$, discarding those that vanish identically. They generate $\C[N]$ and have root characters $\alpha_i$. Let $V$ be their rational span, of dimension $r$. Integrality gives a point at which all these functions are nonzero. Its stabilizer is their common kernel, so its orbit has dimension $r$; in particular $r\le d$.

Let $H$ be the identity component of that kernel. A character is trivial on $H$ precisely when it belongs to $V$ over $\Q$. Thus the roots of $G'=Z_G(H)$ are $\Phi\cap V$, and they span $V$ because they include the $\alpha_i$. Choose a cocharacter $\lambda$ of $H$ avoiding the finitely many root hyperplanes that do not contain $H$. Then $G'=Z_G(\lambda)$ is a Levi subgroup of derived rank $r$. If $H$ is trivial, take $\lambda=0$ and $G'=G$; the same character calculation then gives $r=\rank G$. Since its coordinate generators have weight zero under $\lambda$, $N$ is fixed as a scheme.

Let $S$ be the connected component of $(G/B)^\lambda$ containing $N$. The pattern map gives a $G'$-equivariant isomorphism $\psi:S\longrightarrow G'/B'$, where $B'=B\cap G'$. For every fixed-point label $v$ in this component, it identifies the fixed ordinary cell $(C_v)^\lambda$ and the fixed opposite cell $(C^v)^\lambda$ with $C'_{\phi(v)}$ and $(C')^{\phi(v)}$, respectively \cite[Theorem~10]{BB03}\cite[Theorem~4.1 and the following paragraph]{Woo10}. These are statements about the fixed cells themselves; no irreducibility assertion for $X_w\cap S$ is needed.

The open Richardson stratum $C_w\cap C^u$ is nonempty and dense in $N$. Since it is fixed, the connected fixed cell $(C_w)^\lambda$ meets $S$, and hence lies in $S$; in particular $w\in S$. Likewise $(C^u)^\lambda\subset S$. The two cell isomorphisms are restrictions of the same injective map $\psi$, so they identify their intersection with the entire intersection of their images:
\[
 \psi(C_w\cap C^u)=C'_{\phi(w)}\cap(C')^{\phi(u)}.
\]
Take reduced closures inside the fixed opposite cell on the left and inside $(C')^{\phi(u)}$ on the right. The first closure is $N$, since $N$ is reduced and closed in $C^u$, and the indicated stratum is dense. The second is the Levi Schubert slice $N'_{\phi(u),\phi(w)}$. This gives an isomorphism of actual affine schemes. The same argument applies to every $N_{u,z}\subset N$ with $u\le z\le w$, so the identifications respect closed subslices.

Since $T\subset G'$ and $\psi$ is $G'$-equivariant, it preserves the ambient characters. Conjugation to a standard Levi preserves root lengths. All maps were constructed over $\C$, so the isomorphism persists after any extension of the ground field.
\end{proof}

\subsection{Finite covers and the Schubert boundary}
\begin{lemma}\label{tor:closed-gluing}
Suppose finitely many closed subschemes $T_a\subset T$ cover its support, and $P\longrightarrow T$ is finite \'etale. Sections on the $T_a$ that agree on their schematic pairwise intersections glue uniquely to a section on $T$. It suffices to check agreement on a family covering the support of each intersection.
\end{lemma}
\begin{proof}
The chosen sheet and the union of the other sheets above each $T_a$ are closed in $P$. Agreement makes the finite unions of these two families disjoint, and together they cover $P$. The chosen union is therefore clopen. With its induced scheme structure it is finite \'etale over $T$ and has one point in each geometric fibre, hence is isomorphic to $T$. Finally, the equalizer of two sections of a finite \'etale map is clopen. Agreement on its whole support is therefore sufficient, even if the intersection is nonreduced.
\end{proof}

\begin{theorem}\label{tor:boundary}
Let $d=\ell(w)-\ell(u)\ge4$. A finite \'etale cover of $U_F(u,w)$ that splits on every $U_F(u,v)$ with $u<v\lessdot w$ is split. Moreover, restriction induces an injection
\begin{equation}\label{tor:boundary-injection}
 \Tor\Pic U_F(u,w)\longrightarrow
 \bigoplus_{u<v\lessdot w}\Tor\Pic U_F(u,v).
\end{equation}
\end{theorem}
\begin{proof}
On the simply connected form of the group, a boundary section of $\mathcal L_\rho$ defines an effective Cartier divisor supported on all Schubert divisors, with positive multiplicities \cite{AS14}. The flag variety is unchanged by this choice of isogeny. Restrict the section to $N_{u,w}$; it is nonzero because the slice meets the top cell. Trivializing at the vertex gives a nonzero element $f$ of the maximal ideal of $A=A_F(u,w)$. Write $U=\Spec A\setminus\{s\}$ and $T_0=\Spec(A/fA)\setminus\{s\}$. Its support is covered by the closed coatom pieces $U_v$.

The support of an intersection of lower pieces $U_{z_i}$ is covered by the common lower pieces $U_q$ with $u<q\le z_i$ for every $i$. This follows by restricting the Schubert support intersection identity to the opposite cell, using the common field and the quotient compatibility above. The case $q=u$ disappears on the puncture. No reducedness assertion about these intersections is needed.

The order complex of the open Bruhat interval $(u,w)$ is a sphere $S^{d-2}$ \cite{BW82}; in this range it is connected and simply connected. Because every $U_z$ with $u<z<w$ is connected, a cover of $T_0$ split on each coatom is split on every lower piece. The section sets are finite, and the restriction maps are bijections satisfying the usual composition identities. Transport labels along comparability edges. Backtracking gives the identity, and each triangle gives the identity by composition. Simple connectivity therefore makes transport independent of the path, and all these finite sets can be labelled compatibly.

Sections with the same label agree on every common lower piece, hence on the entire support of each schematic intersection. Lemma~\ref{tor:closed-gluing} glues them on $T_0$, so the restricted cover splits there. The ring $A$ is henselian and has depth $d$. By \cite[Tag~0BLJ]{The26}, restriction from finite \'etale covers of $U$ to those of $T_0$ is fully faithful. Lift a splitting isomorphism and its inverse; faithfulness makes the lifts inverse. This proves the first assertion. Essential surjectivity is not used.

For the second assertion, take $L$ in the kernel and choose $L^{\otimes m}\simeq\OO_U$. Its Kummer torsor restricts to a split torsor on each $U_v$: since $L|_{U_v}$ is trivial, the restricted torsor class can only come from a unit modulo $m$th powers. Normal Hartogs extension identifies these units with $A_F(u,v)^*$, and every such unit has an $m$th root by strict henselianity in characteristic zero. The first assertion splits the original torsor, so $L$ is trivial.
\end{proof}

\subsection{An integral curve-degree bound}
\begin{theorem}\label{tor:degrees}
Let $A$ be a Noetherian strictly henselian normal local $\Q$-algebra of dimension at least two, and put $U=\Spec A\setminus\{s\}$. Let $f:Y\longrightarrow\Spec A$ be proper birational, with $Y$ regular and integral. Suppose every finite \'etale cover of its connected closed fibre $F_0$ splits. Let $E_1,\ldots,E_r$ be all prime divisors of $Y$ contained in $F_0$. For actual proper integral $1$-cycles $Z_i$ on $F_0$, allowing signed coefficients, set
\[
 N_{ij}=\deg(\OO_Y(E_j)|_{Z_i}).
\]
If $N$ has full column rank and largest Smith invariant $h$, then $h\Tor\Pic(U)=0$. For $r=0$, take $h=1$.
\end{theorem}
\begin{proof}
Choose $L^{\otimes m}\simeq\OO_U$. Its generic Kummer algebra is $K[z]/(z^m-a)$ for $a\in K^*$, where $K=\Frac A$; this algebra need not be a field. The pulled-back torsor is \'etale on $Y\setminus F_0$. At every prime divisor not contained in $F_0$, the DVR description of a Kummer torsor shows that the valuation of $a$ is divisible by $m$. Consequently
\[
 \operatorname{div}_Y(a)=mD+\sum_jc_jE_j
\]
for an integral divisor $D$ and integers $c_j$. All these divisors are Cartier because $Y$ is regular.

Taking degrees on the $Z_i$ gives $Nc\equiv0\pmod m$. Smith normal form supplies an integer matrix $M$ with $MN=hI_r$. Hence $hc=mb$ for an integer vector $b$, and
\[
 \operatorname{div}_Y(a^h)=mD',
 \qquad D'=hD+\sum_jb_jE_j.
\]
Locally write $D'=\operatorname{div}(g_\alpha)$. The algebras
\[
 \OO_Y[t_\alpha]/(t_\alpha^m-a^h/g_\alpha^m),
 \qquad t_\alpha=(g_\beta/g_\alpha)t_\beta,
\]
glue by the indicated transitions to a finite \'etale Kummer torsor on $Y$. Proper henselian equivalence \cite[Tag~0A48]{The26} and the hypothesis on $F_0$ split this torsor. Thus $a^h$ is an $m$th power in $K$. The $h$th power of the original torsor is generically split, hence split on the normal scheme $U$. In particular $L^{\otimes h}$ is trivial.
\end{proof}

Only genuine curve degrees and the complete set of columns are required. The rows need not generate all curve degrees. In particular, a full-column minor equal to $\pm1$, or an integer combination of full-column maximal minors equal to $1$, proves vanishing of torsion.

\subsection{Bott--Samelson sources and exact exceptional columns}\label{tor:sources}
Fix a reduced word $Q=(i_1,\ldots,i_l)$ for $w$, with Bott--Samelson resolution $f_Q:Y_Q\longrightarrow X_w$ and canonical boundaries $D_1,\ldots,D_l$. The cell's unipotent section and $B$-equivariance give the compatible product
\[
 f_Q^{-1}(O_u\cap X_w)\simeq C_u\times Y_N,
 \qquad Y_N=f_Q^{-1}(N_{u,w}).
\]
The product is obtained by acting on $Y_N$ with the unique cell coordinate in the chosen unipotent section, and it transports the boundary bundles. Its left side is smooth and integral, so $Y_N$ is also smooth and integral; the map $Y_N\longrightarrow N_{u,w}$ is proper birational.

After extension to $F$ and base change to $A=A_F(u,w)$, the resulting $Y_A$ is regular, since localization and ind-\'etale base change preserve regularity. It is integral as well. Indeed, proper birationality to a normal target gives $f_*\OO=\OO$; flat base change preserves this identity \cite[Tag~02KH]{The26}, so $Y_A$ is connected. Distinct irreducible components of a regular Noetherian scheme are disjoint and open, which proves integrality. Its generic fibre remains a point.

\Needspace{5\baselineskip}
\begin{lemma}\label{tor:fibre}
The closed fibre of $Y_A$ is connected, and every finite \'etale cover of it splits.
\end{lemma}
\begin{proof}
The complex Bott--Samelson fibre $F_u=f_Q^{-1}(u)$ is connected by normality of the target. It is affine-paved: its intersections with the projective Bia\l ynicki--Birula cells are the polynomial graphs of \cite[Proposition~2.1]{Ha04}, and the closed Bia\l ynicki--Birula filtration restricts to the fibre.

For a connected proper complex affine-paved scheme, compactly supported cohomology of the paving has only even degrees. The number of zero-dimensional cells is therefore its zeroth Betti number, namely one. A finite \'etale cover of degree $n$ pulls each affine cell back to $n$ copies, because complex affine space has no nontrivial finite \'etale cover by Riemann existence \cite[Expos\'e~XII, Theorem~5.1]{Gro71}. The cover is proper and has $n$ zero-dimensional cells. It consequently has $n$ connected components, each of degree one, and is split.

Nilpotents do not affect the finite \'etale category. The closed fibre over $A$ is $F_u\times_{\C}F^{\mathrm{sep}}$. In characteristic zero, $F^{\mathrm{sep}}$ is algebraically closed. Proper invariance of the finite \'etale category under extension of algebraically closed fields \cite[Tag~0A49]{The26} proves the assertion.
\end{proof}

\begin{lemma}\label{tor:columns}
Set $z_j=\Dem(Q\setminus j)$. The prime divisors of $Y_A$ contained in its closed fibre are precisely the restrictions of the $D_j$ with $z_j=u$, each reduced and with multiplicity one. In particular the exact column set is $E(Q,u)=\{j:\Dem(Q\setminus j)=u\}$.
\end{lemma}
\begin{proof}
The divisor $D_j$ is the smooth integral Bott--Samelson variety for the deleted word, whether or not that word is reduced; its image is $X_{z_j}$. If $u\le z_j$, its restriction to the source product of Subsection~\ref{tor:sources} has the form $C_u\times D_{j,N}$. If $z_j=u$, the smooth integral $D_{j,N}$ lies entirely over the vertex. It is geometrically integral, so extension to $F$ followed by the common strict base change leaves a single reduced prime divisor. Its bundle is the restriction of $\OO_{Y_Q}(D_j)$.

For completeness, first work over $F$, rather than comparing heights across a possibly transcendental extension of $\C$. A height-one closed-fibre prime on $Y_A$ contracts, by flat going down, to a prime of height at most one on $Y_{N,F}$. It cannot contract to the generic point, whose image is generic in $N_F$ rather than the vertex. It thus contracts to a prime divisor over the vertex. Its product with $C_{u,F}$ closes to a prime divisor on $Y_{Q,F}$ outside the inverse image of the top cell. That complement is the union of the canonical $D_{j,F}$, so this prime is one of them. Its image centre on the slice is $N_{u,z_j}$, which lies over the vertex exactly when $z_j=u$. Distinct columns have distinct contractions, so they cannot merge after strict base change. This proves completeness and the stated multiplicities.
\end{proof}

\begin{lemma}\label{tor:products}
Let $N=N_1\times\cdots\times N_s$, where $s\ge2$ and every $N_i$ is a positive-dimensional actual Schubert slice. Then, for every $F\supset\C$, the strict puncture at the product vertex has torsion-free Picard group. The same conclusion holds if the factors are scalar Veronese cones.
\end{lemma}
\begin{proof}
Take the product of the regular proper birational slice sources just constructed. The closed fibre of each factor is a proper closed subset of an integral source, hence has codimension at least one. Their product has codimension at least two and contains no prime divisor. It is connected and affine-paved, so its finite \'etale covers split by the preceding argument.

After extension to $F$ and strict base change, regularity and integrality persist as above. A new height-one closed-fibre prime would contract by flat going down to a prime of height at most one in the original product source, contradicting the codimension of its vertex fibre. The empty-column case of Theorem~\ref{tor:degrees} applies. For a scalar Veronese cone of dimension $m\ge2$, use the regular source $\operatorname{Tot}(\OO_{\PP^{m-1}}(-2))$; a one-dimensional factor is smooth and uses the identity. The same product argument has projective spaces or points as its closed fibres.
\end{proof}

\subsection{Mask characters as integral curve degrees}
Use the adjoint form for this calculation, so that $X^*(T)$ is the root lattice and every root is primitive. This change leaves the flag variety and the Bott--Samelson variety unchanged; it is independent of the simply connected form used for $\mathcal L_\rho$ above.

For a mask $\gamma\in\{0,1\}^l$, write $\gamma^j=s_{i_1}^{\gamma_1}\cdots s_{i_j}^{\gamma_j}$. Its endpoint is $\gamma^l$. The canonical linearization on $\OO(D_j)$ has fibre character
\begin{equation}\label{tor:characters}
 \chi_j(\gamma)=
 \begin{cases}
  0,&\gamma_j=1,\\
  -\gamma^{j-1}\alpha_{i_j},&\gamma_j=0.
 \end{cases}
\end{equation}
Indeed, the invariant boundary section is nonzero at a selected position; at an omitted position its normal character is the transported negative simple root. The bundle is pulled back through later stages of the tower, so later coordinates do not change this character.

\begin{proposition}\label{tor:mask-degrees}
Fix an endpoint-$u$ reference mask $\gamma_0$ and simple-root coordinate functionals $\lambda_a$ on the root lattice. Every row
\begin{equation}\label{tor:raw-rows}
 \bigl(\lambda_a(\chi_j(\gamma)-\chi_j(\gamma_0))\bigr)_{j\in E(Q,u)},
 \qquad \gamma^l=u,
\end{equation}
is the simultaneous degree row of a signed integral proper $1$-cycle on the closed fibre. Consequently, a full-column-rank matrix of these rows with largest Smith invariant $h$ annihilates $\Tor\Pic U_F(u,w)$ for every $F\supset\C$.
\end{proposition}
\begin{proof}
First join endpoint-$u$ masks by moves changing two bits and keeping the endpoint. Coxeter deletion successively removes two selected letters from a nonreduced selected expression, ending at a reduced $u$-subword. These reduced subwords are complements of the facets of the subword complex $\Delta(Q,u)$. Its pure shellability \cite[Theorems~2.5 and 3.7]{KM04} makes its facet adjacency graph connected; an adjacency exchanges two positions. Thus the entire endpoint-mask graph is connected.

For one such move, let $p<q$ be its changed positions, and let $\beta$ be the positive representative of $\pm\gamma^{p-1}\alpha_{i_p}$. Direct multiplication gives
\[
 \delta^t=
 \begin{cases}
  \gamma^t,&t<p\text{ or }t\ge q,\\
  s_\beta\gamma^t,&p\le t<q.
 \end{cases}
\]
Equality of the endpoints forces $\gamma^{q-1}s_{i_q}(\gamma^{q-1})^{-1}=s_\beta$. Thus both changed positions are $\beta$-walls. The masks agree at every other position and are equivalent for the codimension-one subtorus $\ker\beta$ in the terminology of \cite[Sections~1 and 4]{Ha04}.

Here is the precise fibre geometry of this equivalence class. List all its $\beta$-walls as $a_1<\cdots<a_m$, and let $G_\beta$ be the rank-one subgroup generated by $T$ and the root subgroups for $\pm\beta$. After choosing the minimal gallery in the class, the partial-gallery parametrization of \cite[equation~(4.1), Proposition~4.2 and Corollary~4.4]{Ha04} identifies its fixed component with $(G_\beta/B_\beta)^m\simeq(\PP^1)^m$. Explicitly, between walls $a_j$ and $a_{j+1}$, the partial flag is $g_j\gamma_{min}^tB$, with $g_jB_\beta$ the $j$th coordinate. The flags at the walls recover these coordinates, giving an isomorphism of schemes, not merely a parametrization of closed points.

The final flag is $g_m\gamma_{min}^lB$. Hence the endpoint map is the last projection followed by the isomorphism onto the rank-one orbit in $G/B$. Its fibre over $u$ is exactly $(\PP^1)^{m-1}\times\{u\}$. In particular $m\ge2$ for the two distinct masks under consideration. Successively changing the first $m-1$ coordinates joins these masks by coordinate-axis projective lines in the actual fibre $f_Q^{-1}(u)$.

On each such line, the tangent character at a fixed endpoint is $\beta$ or $-\beta$. For a linearized line bundle, its two fibre characters differ by its degree times the negative of the tangent character at the starting endpoint. This follows from the transition function for $\OO_{\PP^1}(n)$ on the two standard affine charts. Thus, for an oriented axis path $C_1,\ldots,C_t$, there are signs $\eps_\nu\in\{1,-1\}$, independent of the line bundle, such that the signed integral cycle $Z_e=\sum_\nu\eps_\nu[C_\nu]$ satisfies
\[
 \chi_j(\delta)-\chi_j(\gamma)
 =\beta\deg(\OO(D_j)|_{Z_e})
\]
simultaneously for all $j$. Along a path from $\gamma_0$ to $\gamma$, apply $\lambda_a$ and add. This yields \eqref{tor:raw-rows}, realized by the cycle $\sum_e\lambda_a(\beta_e)Z_e$. Every coefficient is integral. The cycles lie over $u$, hence in the slice source, and the product chart transports their boundary bundles. Degrees are unchanged after extension to the algebraically closed residue field. Apply Theorem~\ref{tor:degrees} and Lemmas~\ref{tor:fibre}--\ref{tor:columns}.
\end{proof}

The full raw row lattice also equals the lattice of the divided two-bit-move rows $((\chi_j(\delta)-\chi_j(\gamma))/\beta)_j$. One inclusion follows from the path argument. Conversely, a primitive root $\beta$ has relatively prime simple-root coordinates, so an integral B\'ezout combination of its coordinate differences recovers the divided row. Edge differences are differences of reference rows. This equality survives projection to a fixed column set. Rational normalization of rows, or replacement by a different ambient lattice, would not preserve the argument.

\paragraph{A small degree witness.}
In type $A_2$, take $Q=(1,2,1)$, $u=s_1$, and endpoint masks $\gamma=100$ and $\delta=001$. The full exceptional column set is $E(Q,u)=\{2\}$. Formula~\eqref{tor:characters} gives $\chi_2(\gamma)=-(\alpha_1+\alpha_2)$ and $\chi_2(\delta)=-\alpha_2$, so their difference is $\alpha_1$. The two $\alpha_1$-walls are positions $1$ and $3$; their fixed component is $\PP^1\times\PP^1$, and its endpoint-$u$ fibre is a coordinate line. With the sign prescribed above, this line realizes the raw simple-root rows $(1)$ and $(0)$, the latter by the zero cycle. Thus the undivided character coordinates already give integral degree rows.

\subsection{Passage to ordinary points and partial flags}
We use flat compatibility of duals and reflexive powers from Lemma~2.2.

\begin{lemma}\label{tor:strict-injection}
For an excellent normal local domain $A$ of dimension at least two, pullback to $A^{\sh}$ is injective on the Picard group of the punctured spectrum.
\end{lemma}
\begin{proof}
The closed fibre of $A\longrightarrow A^{\sh}$ is a field, so the upper puncture is the inverse image of the lower puncture. Extend a line bundle on the latter to its unique finite rank-one reflexive module $M$ on $\Spec A$: take a rational section and close its finite Weil divisor. Normality and the codimension of the complement give uniqueness. By Lemma~2.2, its pullback is reflexive. If the pulled-back line bundle is trivial, normal Hartogs uniqueness makes this module $A^{\sh}$. Faithfully flat descent makes $M$ free of rank one, proving injectivity.
\end{proof}

\begin{proposition}\label{tor:ordinary}
For a full-flag Schubert variety $X$, ordinary torsion detection at all scheme points is equivalent to vanishing of punctured Picard torsion at $\eta_Y$ for every Schubert subvariety $Y\subseteq\Sing(X)$ that is not an irreducible component of $\Sing(X)$. It suffices to prove this vanishing after strict henselization.
\end{proposition}
\begin{proof}
First suppose these punctured groups vanish. A local class at $x$ extends to a global Weil class by taking closures of its prime divisors. The affine big cell has trivial class group, so divisor localization represents the global class by a $B$-stable boundary divisor $D$. Suppose its class at $x$ is torsion, all maximal singular generic localizations through $x$ vanish, but $D$ is not Cartier at $x$.

Its closed non-Cartier locus is $B$-stable. A component $Y$ containing $x$ is therefore a Schubert subvariety. At $\eta_Y$, the class is nonzero torsion and is Cartier on the entire puncture of $\OO_{X,\eta_Y}$: a nonfree point there would give a strictly larger nonfree subvariety containing $Y$. Regular local rings are factorial, so $Y\subseteq\Sing(X)$ and $\codim_XY\ge2$. It is not a maximal singular component, because the corresponding detector was assumed zero. Reflexive Hartogs extension shows that triviality on this puncture would imply Cartierness at $\eta_Y$. Thus we have obtained a nontrivial torsion line bundle on a puncture whose Picard torsion vanishes by hypothesis, a contradiction.

Conversely, extend a torsion line bundle on the puncture at such an $\eta_Y$ reflexively. A reflexive power of this extension is trivial by the same extension argument. Every maximal singular component containing $Y$ strictly contains it, so its generic point belongs to the puncture and has zero detector. Torsion detection at $\eta_Y$ makes the extension free, proving the required vanishing. The final assertion is Lemma~\ref{tor:strict-injection}.
\end{proof}

\begin{proposition}\label{tor:partial}
If ordinary torsion detection holds at all scheme points of every full-flag Schubert variety, it also holds for every partial-flag Schubert variety.
\end{proposition}
\begin{proof}
For $X\subset G/P$, its inverse image $Y\subset G/B$ is a Schubert variety, and $\pi:Y\longrightarrow X$ is smooth, surjective and Zariski locally trivial, with geometrically irreducible fibre $P/B$. Thus $\Sing Y=\pi^{-1}(\Sing X)$, and maximal singular components correspond by inverse image. Indeed, the inverse image of an irreducible base under this locally trivial map has irreducible fibres and is irreducible; no vertical component occurs.

At an arbitrary point $x$, choose in a local trivialization the generic point $y$ of an affine big cell of its fibre. Writing $A=\OO_{X,x}$, its ring is the faithfully flat local extension
\[
 B=A[t_1,\ldots,t_a]_{\mathfrak m_AA[t_1,\ldots,t_a]}.
\]
If $\mathfrak p\subset\mathfrak m_A$ is the prime of a maximal singular component through $x$, the corresponding upstairs generic localization is
\[
 B_{\mathfrak pB}
 =A_{\mathfrak p}[t_1,\ldots,t_a]_{\mathfrak pA_{\mathfrak p}[t_1,\ldots,t_a]}.
\]
These are actual compatible detector maps. A finite rank-one reflexive module representing a downstairs torsion class remains reflexive and torsion over $B$ by Lemma~2.2. Vanishing of its downstairs detectors gives vanishing upstairs. Detection on $Y$ makes it free over $B$, and faithfully flat descent makes it free over $A$. No isomorphism of the entire class groups is required.
\end{proof}

\section{The finite-type torsion theorem}\label{tor:finite}

Retain the objects of Section~\ref{tor:local}, and abbreviate $T_F(N)=\Tor\Pic(\Spec R_F(N)^{\sh}\setminus\{s\})$. Every assertion below quantifies over all fields $F\supset\C$. We prove the local vanishing required by Proposition~\ref{tor:ordinary}. The only computational inputs are complete bounded calculations with integer certificates.

\subsection{The finite certificates and what they prove}
Choose one reduced word $Q_w$ for each element of the finite Weyl group under consideration. Testing every reduced expression is unnecessary: each pair needs only one regular source. For $u\le w$, use the full column set $E(Q_w,u)$ from Lemma~\ref{tor:columns}. If it is empty, Proposition~\ref{tor:mask-degrees} already gives $T_F(N_{u,w})=0$. Otherwise the rows are the literal integer rows of \eqref{tor:raw-rows}; no rational row normalization is allowed.

Here is the finite algorithm for the bounded inputs. Starting at the identity, breadth-first search by right multiplication in the simple reflections constructs the group and one shortest word for every element. Closure under every generator proves completeness, and the first word found is reduced. For each word, compute every single deletion by the Demazure rule: multiply by a letter exactly when it increases length. Retain the resulting exact centre $u$ and every position producing it. Evaluate masks by ordinary group multiplication and characters by \eqref{tor:characters}. Every displayed or stored minor uses every column, with distinct valid row indices. Calculate determinants and B\'ezout identities over the integers.

\begin{proposition}\label{tor:certificates}
The following bounded cases have primitive full-column row lattice:
\begin{enumerate}
 \item in $B_3$, every exact-centre pair of dimension three;
 \item in $C_3$, every exact-centre pair of dimension three except $(Q_w,Q_u)=(121321,121)$;
 \item in $F_4$, every exact-centre pair of dimension four.
\end{enumerate}
The exceptional $C_3$ pair has one column and raw difference lattice $2\Z$. In rank two, the only possible nonzero type-$C$ surface input is $(Q_w,Q_u)=(121,1)$, again with bound two.
\end{proposition}
\begin{proof}
For $B_3$ and $C_3$, use signed permutations of three coordinates, with simple roots $\alpha_1=e_1-e_2$, $\alpha_2=e_2-e_3$, and $\alpha_3=e_3$ in $B_3$ or $\alpha_3=2e_3$ in $C_3$. Each group has $48$ elements. The algorithm gives $199$ distinct exact-centre pairs in each type: $138$ have gap one, $49$ have gap two, and $12$ have gap three. Thus there are $61$ exceptional exact-centre pairs of gap at least two. Table~\ref{tor:bc-table} lists all twelve gap-three pairs. Each matrix uses all exact columns and has determinant $\pm1$, except for the indicated $C_3$ pair.

For that exceptional pair, the endpoint masks are $001011$, $011010$, $100011$, $110001$, and $111000$. Their exceptional-column characters, in simple-root coordinates, are respectively
\[
 (0,0,-1),\quad (0,-2,-1),\quad (0,0,-1),\quad
 (-2,-2,-1),\quad (-2,-2,-1).
\]
Hence all raw differences are even, and some coordinate difference is $2$ or $-2$. The complete raw difference lattice is therefore $2\Z$.

The rank-two enumeration gives the surface assertion. In the $C_2$ convention, $\alpha_1=e_1-e_2$ is short and $\alpha_2=2e_2$ is long. For $Q=121$ and $u=s_1$, the exact column is position $2$, and the endpoint masks are $100$ and $001$. The character difference in this column is $2(e_1-e_2)$, whose simple-root coordinate column is $(2,0)^{\mathsf t}$ up to overall sign. There are exactly three nonempty dimension-two exact column sets in $C_2$:
\[
\begin{array}{c c c c}
 Q_w&Q_u&E(Q_w,u)&\text{full-column degree matrix}\\[2pt]
 121&1&\{2\}&(2)\\
 212&2&\{2\}&(1)\\
 1212&12&\{2,3\}&\begin{pmatrix}0&-1\\1&-1\end{pmatrix}
\end{array}
\]
For the last pair, reference mask $0011$ and mask $1100$ give exactly the two displayed root-coordinate rows. This is also a direct rank-two check, independent of the enumeration. Every other pair of the stated dimension has empty exact column set. The group operations in $B_3$ and $C_3$ agree, but the character matrices differ because the root lengths differ.

For $F_4$, use simple roots $e_2-e_3$, $e_3-e_4$, $e_4$, and $\tfrac12(e_1-e_2-e_3-e_4)$. The integral reflection matrices in the simple-root basis generate $1152$ elements, with maximum length $24$. Every simple reflection acts within the closed enumerated set, so this is exactly the generated Weyl group. All single deletions give $567$ dimension-four exact-centre pairs: $497$ have one column, $49$ have two, and $21$ have three. The accompanying reconstructor independently generates all these pairs and their raw mask-coordinate rows. For each pair it records full-column minors and integer B\'ezout coefficients giving a combination equal to $1$.

For a dimension-four pair, an endpoint mask selects at least $\ell(w)-4$ letters and has the same parity as $\ell(w)-4$. The complete reduced word cannot have the shorter endpoint. Thus it suffices to enumerate masks omitting two or four letters. Each stored witness is also checked directly by its endpoint and characters. The supplementary calculation retains complete group and centre coverage, all exceptional columns, legal row witnesses, and exact determinant identities. A list of numerical gcds alone would not establish these facts. Proposition~\ref{tor:mask-degrees} now gives the stated consequences for the actual punctured local objects over every $F$.
\end{proof}

\begin{table}[tbp]
\centering
\small
\renewcommand{\arraystretch}{1.3}
\begin{tabular}{@{}ccccc@{}}
\toprule
$Q_w$&$Q_u$&Exact positions&$B_3$ square matrix&$C_3$ square matrix\\
\midrule
$13213$&$13$&$3$&$(-1)$&$(-1)$\\
$121321$&$121$&$4$&$(-1)$&$(-2)$\\
$132132$&$132$&$3$&$(-1)$&$(-1)$\\
$213213$&$213$&$4$&$(1)$&$(1)$\\
$1213213$&$1213$&$4,5$&$\begin{pmatrix}-1&0\\-1&1\end{pmatrix}$&$\begin{pmatrix}-2&1\\1&-1\end{pmatrix}$\\
$1321323$&$1323$&$3$&$(-1)$&$(-1)$\\
$2132132$&$2132$&$4$&$(1)$&$(1)$\\
$2321323$&$2323$&$4$&$(-1)$&$(-1)$\\
$12132132$&$12132$&$4,5$&$\begin{pmatrix}-1&0\\-1&1\end{pmatrix}$&$\begin{pmatrix}-2&1\\1&-1\end{pmatrix}$\\
$12321323$&$12323$&$5$&$(-1)$&$(-1)$\\
$21321323$&$21323$&$4$&$(1)$&$(1)$\\
$121321323$&$121323$&$4,5,6$&$\begin{pmatrix}-1&0&1\\0&0&-1\\-1&1&0\end{pmatrix}$&$\begin{pmatrix}-2&0&1\\1&0&-1\\-2&1&0\end{pmatrix}$\\
\bottomrule
\end{tabular}
\caption{Complete dimension-three exact-centre list in $B_3$ and $C_3$, with witnesses reconstructed from raw character rows. Every displayed matrix has determinant $\pm1$ except $(-2)$. Positions are one-based and comprise the entire exact column set. A digit $1$ in a word denotes a simple reflection, not the identity.}\label{tor:bc-table}
\end{table}

\Needspace{8\baselineskip}
\subsection{Simply laced types and type \texorpdfstring{$B$}{B}}
\begin{proposition}\label{tor:simply-laced-b}
For a simply laced slice, $T_F(N)=0$ in every dimension at least two. For a type-$B$ slice it vanishes in every dimension at least three.
\end{proposition}
\begin{proof}
In dimensions at most three, Proposition~\ref{tor:levi} reduces to a Levi subgroup of rank at most three. If at least two positive-dimensional factors occur, apply Lemma~\ref{tor:products}. The irreducible simply laced factors are $A_1$, $A_2$, and $A_3$. The first two have only smooth slices. In $A_3$, the only singular three-dimensional pairs are $(u,w)=(1324,3412)$ and $(2143,4231)$. Their rank equations reduce respectively to $c=0$, $ae+bd=0$, and to $bc-ad=0$, up to independent smooth variables; at dimension three these are ordinary threefold nodes.

A regular small source of the node is $\operatorname{Tot}(\OO_{\PP^1}(-1)\oplus\OO_{\PP^1}(-1))$. Its vertex fibre is $\PP^1$ and contains no divisor, so the empty-column case of Theorem~\ref{tor:degrees} kills torsion. These low-rank equations can also be read from the standard generic-singularity descriptions \cite{Man01a,Cor03,BW03}. As an independent check, the supplementary reconstructor verifies all dimension-two and dimension-three exact-centre inputs in $A_3$; their row lattices are primitive, and all omitted pairs have empty column set.

For a single type-$B$ factor, the rank-three input of Proposition~\ref{tor:certificates} gives the dimension-three base. Smaller standard Levi factors are included: restricting the signed-permutation computation to a smaller alphabet gives the complete rank-one or rank-two computation. Empty column sets are handled by Proposition~\ref{tor:mask-degrees}.

For $d\ge4$, Theorem~\ref{tor:boundary} injects the dimension-$d$ group into the dimension-$(d-1)$ coatom groups. Induction proves both assertions. Every low-dimensional model and source was constructed over $\C$ and transported to the same $F$, so the induction never compares different cell function fields.
\end{proof}

\subsection{Type \texorpdfstring{$C$}{C}: the equivariant induction statement}
The remaining classical case needs more than local vanishing. Use signed ambient coordinates $\eps_i$ for the type-$C$ roots. A half-character $\eps_i$ need only belong to a finite covering torus; the actual coordinate characters below are roots $2\eps_i$ and $\eps_i+\eps_j$.

\begin{theorem}\label{tor:veronese-classification}
Let $N$ be an actual complex type-$C$ Schubert slice of dimension $d$. If $T_F(N)\ne0$ for some $F\supset\C$, then $d\ge2$ and $N$ is isomorphic, preserving its vertex and ambient torus action, to the scalar Veronese cone $V_d=\A^d/\{\pm1\}$, where the involution acts by simultaneous negation. Its coordinate axes have characters given by $d$ distinct ambient long roots. For this model, $T_F(N)=\Z/2$ for every $F\supset\C$.
\end{theorem}

The quantifier over every $F$ and the equivariance are part of the induction statement. An abstract analytic model, or a formal model without its characters, would not suffice in the next step.

\begin{lemma}\label{tor:veronese-base}
Theorem~\ref{tor:veronese-classification} holds for $d\le3$.
\end{lemma}
\begin{proof}
For $d\le1$, the puncture is empty or the spectrum of a field. For $d=2,3$, use the actual Levi reduction and Lemma~\ref{tor:products}. The single type-$A$ factors are covered by Proposition~\ref{tor:simply-laced-b}. The finite inputs leave only the $C_2$ pair $(s_1,s_1s_2s_1)$ and the $C_3$ pair with words $(Q_u,Q_w)=(121,121321)$.

Here is an actual equivariant model for both. In a symplectic space $E\oplus E^*$ with $\dim E=m$, the graph chart of $\LG(m,2m)$ consists of symmetric maps $S:E\longrightarrow E^*$. The Schubert incidence $Z_m=\{L:\dim(L\cap E)\ge m-1\}$ is the reduced rank-at-most-one symmetric-matrix locus in this chart. Choosing a hyperplane in $E$ and then its isotropic extension gives an irreducible $\PP^1$-bundle incidence, so $Z_m$ is integral of dimension $m$. Its top minimal word is $v=s_1\cdots s_m$. The full-flag inverse image has top $w=vw_J$ and lower endpoint $u=w_J$, where $J=\{1,\ldots,m-1\}$. On the opposite cell at $u$, the projection identifies the full-flag slice with the graph chart.

The finite map $x\longmapsto xx^{\mathsf t}$ has reduced image exactly the symmetric matrices of rank at most one. Its image ring is $\C[x_ix_j]$, the even invariant ring. For $m=2,3$, this gives the stated slices, since $123121=121321$. Their matrix characters are $\eps_i+\eps_j$, and their axis characters are $2\eps_i$. Conjugating a Levi changes these only by signed permutation and preserves root lengths. The next lemma supplies the torsion calculation.
\end{proof}

\subsection{A cover that cannot split on a normal subslice}
\begin{lemma}\label{tor:veronese-cover}
For $d\ge2$, the strict puncture of $V_d$ has \'etale fundamental group $\Z/2$ and Picard torsion $\Z/2$, over every $F\supset\C$. Its nontrivial double cover remains nonsplit on the puncture of every positive-dimensional normal local quotient of its strict vertex ring.
\end{lemma}
\begin{proof}
Let $A=A_F(V_d)$, and base change the finite scalar cover $F[x_1,\ldots,x_d]_{(x)}$ to $A$, obtaining $B$. Its closed fibre has one underlying point. Thus $B$ is local; it is strictly henselian and regular by finite henselian decomposition and ind-\'etale preservation of regularity. Before strict henselization, the invariant maximal ideal is generated by the quadratic monomials $x_ix_j$. After base change, $B/(x_1,\ldots,x_d)B$ is the separably closed residue field of $A$. Consequently the maximal ideals satisfy
\[
 \mathfrak n_B=(x_1,\ldots,x_d)B,
 \qquad \mathfrak m_AB=\mathfrak n_B^2,
\]
where the second identity uses the same quadratic generators.

The puncture of $B$ is a connected \'etale double cover of the puncture of $A$. Every finite \'etale cover of the regular strict puncture of $B$ extends to $\Spec B$ by purity \cite[Tag~0BMA, assertion~(1)]{The26}, and splits there. Hence the downstairs fundamental group is $\Z/2$. Normal Hartogs extension and strict henselian roots of units identify $H^1(U,\mu_m)$ with $\Pic(U)[m]$. Indeed, $\Gamma(U,\OO_U)^*=A^*$, and each unit has an $m$th root, so the unit term in the Kummer exact sequence vanishes. Thus $\Pic(U)[m]=\Hom(\Z/2,\mu_m)$: there is exactly one nonzero torsion class, of order two, corresponding to the scalar double cover.

Now let $A\twoheadrightarrow R$ be a positive-dimensional normal local domain. If the cover splits on its puncture, choose a generic branch to obtain $B\otimes_A R\longrightarrow\Frac R$. Integrality and normality put its image in $R$, giving an $A$-algebra map $B\longrightarrow R$. This map is local: the inverse image of $\mathfrak m_R$ lies over $\mathfrak m_A$, and $B$ has a unique prime over $\mathfrak m_A$, namely $\mathfrak n_B$. The quadratic-ideal identity and surjectivity of $A\longrightarrow R$ now imply $\mathfrak m_R=\mathfrak m_AR\subset\mathfrak m_R^2$. Nakayama's lemma contradicts the positive dimension of $R$. This also covers dimension one: the argument propagates a nonsplit cover, not nonzero Picard torsion on a curve.
\end{proof}

\subsection{Propagation along coatoms}
\begin{lemma}\label{tor:coatom}
Assume Theorem~\ref{tor:veronese-classification} in every dimension less than $d\ge4$. If a dimension-$d$ slice has nonzero $T_F$, every coatom slice is an equivariant scalar Veronese cone, and every proper upward root curve at its vertex has long-root character.
\end{lemma}
\begin{proof}
Theorem~\ref{tor:boundary} and induction inject $T_F$ into a sum of groups each equal to zero or $\Z/2$. Thus a nonzero class has order two. Its nontrivial Kummer double cover is nontrivial on at least one coatom. By induction that coatom is an equivariant Veronese cone, and by Lemma~\ref{tor:veronese-cover} the restricted cover is its scalar cover.

Suppose another coatom shares a lower label $z>u$ with this one. The actual shared-field quotient $A_F(u,z)$ is positive-dimensional and normal, so the restricted cover is nonsplit there by Lemma~\ref{tor:veronese-cover}. If its Picard class on the second coatom were zero, normal Hartogs extension and strict roots of units would make the cover split on the whole second-coatom puncture, a contradiction.

The graph of coatoms joined when they share a lower label $z>u$ is connected. Indeed, the comparability graph of $(u,w)$ is connected; choose a coatom above each vertex. The coatoms chosen above comparable vertices share the smaller vertex. Nontriviality therefore propagates to every coatom, and induction identifies them all equivariantly as Veronese cones.

Every proper upward root curve lies in one of these coatom slices. Its effective torus is the full diagonal torus of the scalar Veronese cone, whose one-dimensional orbit closures through the vertex are its coordinate axes. These axes have long-root characters. This argument makes no assumption about a possible direct edge from $u$ to $w$.
\end{proof}

\begin{lemma}\label{tor:short-edge}
In the situation of Lemma~\ref{tor:coatom}, a direct root edge from $u$ to $w$ cannot be short.
\end{lemma}
\begin{proof}
A reflection changes length by an odd integer, so if its gap is $d\ge4$, then $d\ge5$. The short positive roots in $C_n$ are the $D_n$ roots $e_i\pm e_j$ for $i<j$. Let $\ell_D^*(x)$ count the roots in this set inverted by a signed permutation $x$. Then $\ell_C(x)=\ell_D^*(x)+\operatorname{neg}(x)$, where $\operatorname{neg}(x)$ counts negative entries. A short reflection changes this count by at most two, and therefore $\ell_D^*(w)-\ell_D^*(u)\ge3$.

The endpoints lie in the same signed-parity coset. For the even coset, use $W(D_n)$ itself; for the odd coset, multiply on the left by the sign change $\tau$ in the last coordinate. This $\tau$ interchanges $e_i-e_n$ and $e_i+e_n$ and fixes the other positive $D_n$ roots. It preserves that positive system and $\ell_D^*$. For a positive root $\alpha$, the reflection criterion for Bruhat order says $x<s_\alpha x$ if and only if $x^{-1}\alpha$ is positive \cite[Chapter~2]{BB05}. The positive $D_n$ roots are precisely the positive short $C_n$ roots, and signed permutations preserve the set of short roots. Thus reflection-ascent directions agree in the two systems. In the odd coset, conjugate the reflection by $\tau$; the same criterion proves agreement after transport.

The transported endpoints are consequently comparable in $D_n$, with length gap at least three. Choose the first cover of a saturated $D_n$ chain between them. Its upper endpoint is strictly between the two endpoints. Every step transports back to a short-root reflection ascent in $C_n$, so this first endpoint gives a proper upward short-root curve between $u$ and $w$, contrary to Lemma~\ref{tor:coatom}. In ranks two and three, interpret $D_2=A_1\times A_1$ and $D_3=A_3$; rank one has no short roots. Only this order-preserving direction is used, not equality of induced Bruhat orders.
\end{proof}

\subsection{The actual toric ring}
\begin{proof}[Proof of Theorem~\ref{tor:veronese-classification}]
Lemma~\ref{tor:veronese-base} gives the low-dimensional bases. Proceed by dimension, simultaneously over all ranks and fields. For $d\ge4$, a nonzero class makes every upward curve long by Lemmas~\ref{tor:coatom}--\ref{tor:short-edge}. If there are $b$ such curves, Deodhar's upward inequality gives $b\ge d$ \cite[Proposition~2.3]{CK03}: there are at least $\ell(w)-\ell(u)$ reflections $s_\alpha$ with $u<s_\alpha u\le w$. The corresponding projective root curves lie in $X_w$, and their affine pieces through $u$ in the opposite cell are actual root lines in $N$.

The coordinate characters of this attractive chart lie in a pointed root cone, so both signs of a long axis cannot occur. Distinct long-root curves therefore have independent signed axis directions. Their characters belong to the effective weight span of $\C[N]$, giving $d\le b\le\rank T_{\mathrm{eff}}\le d$.
For the last inequality, choose a point where every nonzero root-coordinate generator is nonzero; its stabilizer is the action kernel. Equality proves that $T_{\mathrm{eff}}$ has a dense open orbit. Since $N$ is normal, it is an actual normal affine toric variety, not merely a toric tangent cone. We use the standard semigroup and orbit--face descriptions \cite[Chapters~1 and 3]{CLS11}.

Discard opposite-cell root coordinates that vanish identically on $N$, and normalize each remaining coordinate to have value one at a complex point of the dense orbit. On that orbit, each coordinate is its weight character. Restriction to the orbit is injective, and these coordinates generate $\C[N]$. Thus $\C[N]=\C[S]$, where $S$ is the semigroup generated by the actual root characters, and $M=\operatorname{gp}S$ is its actual character lattice.

Let $\sigma^\vee=\mathbf R_{\ge0}S$ be the character cone. For a face $\theta$ of this cone, quotienting by all monomials whose weights lie outside $\theta$ gives $\C[S\cap\theta]$. Consequently its one-dimensional faces correspond precisely to invariant irreducible curves through the vertex. Viewed in $G/B$, these are the upward root curves already counted. Their rays are exactly the $d$ independent long axes given by the preceding dimension equality, so the character cone is their positive orthant. The only ambient root-coordinate characters in that orthant are $2\eps_i$ and $\eps_i+\eps_j$.

Each $2\eps_i$ occurs, by restricting its root coordinate to the actual root line. Join $i$ to $j$ if a cross weight $\eps_i+\eps_j$ occurs, and let $J_a$ be the connected components of this graph. The lattice contains $2\Z^d$. Modulo that sublattice, edge vectors span the even-total subspace on each connected component, as follows by using paths in a spanning tree. Therefore
\[
 M=\left\{b\in\Z^d:
   \sum_{i\in J_a}b_i\equiv0\pmod2
   \text{ for every }a\right\}.
\]
Normality of $\C[S]$ implies $S=M\cap\Z_{\ge0}^d$. More explicitly, every $b$ in this intersection has a positive multiple in $S$: because all $2\eps_i$ occur, already $2b\in S$. The character of $b$ is integral over $\C[S]$ in its fraction field, hence belongs to $\C[S]$. This uses normality in $M$, without replacing $M$ by a larger ambient saturation.

This description of $M$ identifies the actual algebra, with its characters, as the tensor product of scalar Veronese invariant rings on the $J_a$. If two or more factors occur, Lemma~\ref{tor:products} gives $T_F(N)=0$, a contradiction. There is one factor, so $N=V_d$ equivariantly. Lemma~\ref{tor:veronese-cover} computes its group over every $F$ and completes the induction.
\end{proof}

\paragraph{A four-dimensional model.}
The rank-at-most-one symmetric $4\times4$ matrices form $\Spec\C[x_ix_j:1\le i\le j\le4]=V_4$. Its character lattice is $\{b\in\Z^4:\sum_i b_i\text{ is even}\}$, whereas its character cone is the full positive orthant. The four axis coordinates have weights $2\eps_i$, and the six off-diagonal coordinates have weights $\eps_i+\eps_j$. Setting all coordinates involving a fixed index equal to zero gives a closed $V_3$. Its scalar double cover remains nonsplit on each such punctured face, as required by the propagation argument. The distinction between the even-sum lattice and $\Z^4$ is essential: saturation in $\Z^4$ would give $\A^4$ and erase the order-two obstruction.

\begin{corollary}\label{tor:c-nonmaximal}
Let $X$ be a type-$C$ full-flag Schubert variety and let $Y\subseteq\Sing(X)$ be a Schubert subvariety that is not an irreducible component of $\Sing(X)$. The strict punctured Picard group of $\OO_{X,\eta_Y}$ has no torsion.
\end{corollary}
\begin{proof}
The only possible nonzero case in Theorem~\ref{tor:veronese-classification} is a slice isolated at its vertex. By the actual product chart of Subsection~\ref{tor:charts}, its label would then be the generic point of a maximal singular component. It cannot be a nonmaximal singular label.
\end{proof}

\Needspace{5\baselineskip}
\subsection{Exceptional types and the conclusion}
\begin{proposition}\label{tor:exceptional}
In $F_4$, $T_F(N)=0$ whenever $\dim N\ge4$. In $G_2$, it vanishes whenever $\dim N\ge3$. Every remaining nonzero group needed here occurs only at the generic point of a maximal singular component.
\end{proposition}
\begin{proof}
For $F_4$ in dimensions at most three, the actual Levi factors are of type $A$, $B$, $C$, or products. The preceding results leave only isolated scalar Veronese models as possible nonzero groups. Here a two-dimensional type-$B$ slice reduces again to rank at most two; the rank-two root system $B_2$ is isomorphic to $C_2$ after interchanging the node labels, and is covered by the surface calculation and graph model above. This also treats rank-two Levi factors inside $B_3$. Dimension four is Proposition~\ref{tor:certificates}, including the empty-column cases. Higher dimensions follow inductively from Theorem~\ref{tor:boundary}.

In $G_2$, reduced words alternate and have length at most six. Deleting an end letter gives length one less. Deleting an interior letter creates exactly one repeated adjacent letter; Demazure idempotence removes the repetition and leaves an alternating word of length two less. Thus single-deletion centres have gap one or two, and those of codimension at least two have gap exactly two. For gap at least three the column set is empty, so Proposition~\ref{tor:mask-degrees} gives vanishing. Normality implies that a nonmaximal singular label has local dimension at least three: in dimension two, normality leaves only an isolated singular vertex in the attractive slice. This covers every required $G_2$ label.
\end{proof}

\begin{theorem}\label{tor:torsion}
For every finite crystallographic root system and every Schubert variety $X\subset G/P$, the ordinary localization map
\[
 \Cl(\OO_{X,x})\longrightarrow
 \bigoplus_{\substack{Z\in\Irr(\Sing X)\\x\in Z}}
 \Cl(\OO_{X,\eta_Z})
\]
is injective on the torsion subgroup, for every scheme point $x$.
\end{theorem}
\begin{proof}
For an irreducible root system, the preceding propositions and Corollary~\ref{tor:c-nonmaximal} kill strict punctured Picard torsion at every nonmaximal singular label. For a reducible root system, the actual slices are products. At least two active positive-dimensional factors are treated by Lemma~\ref{tor:products}; a single active factor is already covered. Take $F=\C(C_u)$ in the generic-point description of Subsection~\ref{tor:charts}, and apply Proposition~\ref{tor:ordinary} to obtain detection at every ordinary point of a full-flag Schubert variety. Proposition~\ref{tor:partial} descends this assertion to all partial flags. This passage uses neither invariance of the entire class group nor any classification in a parabolic Bruhat quotient interval.
\end{proof}

Theorem~\ref{tor:torsion} proves Theorem~1.3.

\section{Cartier indices and canonical loci}\label{sec:consequences}

The two detection theorems now give the index and locus statements. These deductions also make their different scopes explicit.

\begin{proof}[Proof of Corollary~\ref{cor:indices}]
Set $\alpha=[D]\in\Cl(\OO_{X,x})$. The hypothesis that $D$ is $\Q$-Cartier at $x$ means that $\alpha$ has finite order, equal to $\ind_x(D)$. Theorem~\ref{thm:maintorsion} makes its order equal to the order of $\rho_x(\alpha)$. In a finite direct sum the order of an element is the least common multiple of its component orders. This proves the index formula, including the empty case.

Schubert varieties are Cohen--Macaulay, so the canonical sheaf is invertible precisely at the Gorenstein points \cite[Chapter~3]{BH93}\cite[Tag 0DW9]{The26}. Apply the index formula to $K_X$ at a $\Q$-Gorenstein point. Both $\NonQGor(X)$ and $\Bbad(X)$ are contained in $\NGor(X)$. Outside their union the local canonical class is torsion and all its detector images vanish; torsion detection makes it zero. This proves the asserted equality of canonical loci.
\end{proof}

\begin{corollary}\label{cor:typeA-locus}
For type~A Schubert varieties, $\NGor(X)=\Bbad(X)$ without a $\Q$-Gorenstein hypothesis. More generally, \eqref{eq:noncartier} holds for every Weil divisor.
\end{corollary}
\begin{proof}
If $D$ is non-Cartier at $\eta_Z$, it is non-Cartier at every specialization in $Z$, because the free locus is open. Conversely, if $D$ is Cartier at all maximal singular generic points through $x$, Theorem~\ref{thm:mainA} makes its local class at $x$ zero. Taking $D=K_X$ proves the canonical assertion.
\end{proof}

\begin{corollary}[The factorial locus in type A]\label{cor:factorial}
For a type~A Schubert variety,
\[
 \{x:\OO_{X,x}\text{ is not factorial}\}
 =\bigcup_{\substack{Z\in\M(X)\\\OO_{X,\eta_Z}\text{ is not factorial}}}Z.
\]
\end{corollary}
\begin{proof}
If all target class groups through $x$ vanish, injectivity of $\rho_x$ makes $\Cl(\OO_{X,x})$ vanish. A normal Noetherian local domain is factorial exactly when its class group is zero \cite{Fos73}. Conversely, if $\OO_{X,x}$ is factorial, each localization $\OO_{X,\eta_Z}$ is factorial.
\end{proof}

This is the geometric locus conclusion in \cite[Conjecture~6.14]{WY08}\cite[Conjecture~7.39]{WY23}. The global permutation-avoidance classification in \cite{BMB07} addresses a different assertion.

\begin{remark}\label{rem:source-torsion}
The index formula assumes that the local class $[D]\in\Cl(\OO_{X,x})$ is torsion. It neither tests $\Q$-Cartierness from the target indices nor assigns a finite index to an infinite-order class. In particular, generic indices equal to one do not remove $\NonQGor(X)$ from the identity in Corollary~\ref{cor:indices}. The next section exhibits that obstruction.
\end{remark}

\section{A family of counterexamples in type D}\label{sec:type-d}

We prove Theorem~\ref{thm:mainD} by an incidence construction independent of torsion detection and all finite certificates. Fix $r\ge k\ge3$ and put $p=r-k+2$. Let $V$ be a split quadratic space of dimension $2r+4$ with isotropic flag $E_\bullet$. Set $R=E_r$, $F_0=R^\perp$, $P=E_p$, and $H_0=F_0/R$. The quadratic space $H_0$ has dimension four. Choose the ruling of maximal isotropic planes containing $E_{r+2}/R$ and denote its parameter space by $B_0\simeq\PP^1$.

Identify $H_0$ with $U_2\otimes W_2$, with determinant quadratic form. The planes of the chosen ruling are $A_x=\C x\otimes W_2$ and their universal bundle is $\mathcal A\simeq\OO_{\PP^1}(-1)\otimes W_2$. Let $\mathcal E$ be its inverse image in $F_0\otimes\OO_{B_0}$. Choosing a splitting gives $\mathcal E\simeq(R\otimes\OO)\oplus\mathcal A$, of rank $r+2$ and determinant $\OO(-2)$. Every fibre $\mathcal E_A$ is isotropic in $V$.

\subsection{The incidence source and its Schubert image}
On $Z=\Gr_{B_0}(k,\mathcal E)$, let $\mathcal L$ be the tautological subbundle and $\mathcal Q=\mathcal E/\mathcal L$, of rank $p$. Define
\begin{equation}\label{eq:D-incidence}
 \begin{gathered}
 I=\{(A,L)\in Z:\det(P\longrightarrow\mathcal Q)=0\},
 \\ q:I\longrightarrow\OG(k,V),\quad (A,L)\longmapsto L,
 \end{gathered}
\end{equation}
and let $X=q(I)_{\red}$.

\begin{lemma}\label{D:source}
The variety $I$ is integral, normal, and Gorenstein, of dimension $kp$. Its singular locus is the irreducible codimension-three locus $I_2=\{(A,L):\dim(L\cap P)\ge2\}$. Its generic transverse singularity is the node $ad-bc=0$.
\end{lemma}
\begin{proof}
Choose a line $K\subset P$ and then a $k$-plane $L\subset\mathcal E_A$ containing it. This is an integral Grassmann bundle over $\PP(P)\times B_0$, mapping onto the support of $I$. The open set with $\dim(L\cap P)=1$ is nonempty and has a unique such line, so this support is irreducible. A corank-one determinant chart has a nonzero linear parameter, proving generic reducedness. The nonzero determinant section cuts out a Cartier hypersurface in the smooth integral variety $Z$. It is Cohen--Macaulay with no embedded associated primes, so generic reducedness and irreducible support make it integral. Since $\dim Z=kp+1$, its dimension is $kp$.

At a point with $\dim(L\cap P)=c$, the map $P\to\mathcal Q$ has corank $c$. Splitting off an invertible block in an ordinary Grassmann chart leaves the determinant of a matrix of $c^2$ independent parameters; the other coordinates are smooth parameters. One can see the independence directly by decomposing the ambient vector space into $L\cap P$, complements in $L$ and $P$, and a complement of $L+P$. The graph maps from $L\cap P$ to the last complement give precisely that $c\times c$ block. The singular locus is therefore exactly $c\ge2$ \cite{BV88,BH93}.

Choosing a two-plane in $P$ and then a $k$-plane containing it gives an integral incidence variety onto $I_2$. The generic intersection is two-dimensional, where the chosen two-plane is unique, and the dimension is $1+2(p-2)+(k-2)p=kp-3$. Thus $I_2$ is irreducible of codimension three. The hypersurface is regular in codimension one, so Serre's criterion gives normality. It is Gorenstein as a hypersurface. At the generic point of $I_2$, the independent block has size $2\times2$, giving the claimed node.
\end{proof}

\begin{lemma}\label{D:image}
The image $X$ is a Schubert variety, and $q:I\to X$ is proper birational. For $T=\Gr(k,R)$,
\begin{equation}\label{eq:D-small}
 \begin{gathered}
 q^{-1}(T)=T\times\PP^1,\qquad\codim_Iq^{-1}(T)=2k-1,\\
 I\setminus q^{-1}(T)\xrightarrow{\ \sim\ }X\setminus T.
 \end{gathered}
\end{equation}
In particular, $q$ has no exceptional prime divisors.
\end{lemma}
\begin{proof}
The map is proper. Its image is irreducible, closed, and $B$-stable: the chosen Borel preserves $P$, $R$, and $F_0$, and its connectedness makes it preserve the chosen ruling. The orthogonal Grassmannian has finitely many $B$-orbits. An irreducible closed $B$-stable subvariety is therefore the closure of a single orbit, so $X$ is a Schubert variety. In particular, it is normal and Cohen--Macaulay \cite{Ram85,Bri05}.

Let $\pi:F_0\to H_0$ be the quotient. If $L\not\subset R$, its image $\pi L$ is a nonzero isotropic line or plane. It has a unique maximal isotropic extension in the specified ruling. This inverse construction is a morphism: write a nonzero residual vector as a rank-one $2\times2$ matrix and recover its left factor from a nonzero column. These formulae are regular on their nonzero-column opens. For any two residual vectors the left factors agree, because their span belongs to one ruling plane; hence the formulae glue. They give an inverse on the reduced schemes, proving the open isomorphism in \eqref{eq:D-small}.

If $L\subset R$, every ruling plane is allowed, and $\dim(L\cap P)\ge k+p-r=2$. The determinant condition then vanishes identically on the product $T\times B_0$, so the fibre identity is scheme-theoretic. Since $\dim T=k(r-k)$, the exceptional locus has codimension $kp-k(r-k)-1=2k-1\ge5$ in $I$. The open isomorphism proves birationality and absence of exceptional prime divisors.
\end{proof}

For an explicit Weyl-group label, set $n=r+2$ and use positive roots $e_i-e_j,e_i+e_j$ for $i<j$, with simple roots $\alpha_i=e_i-e_{i+1}$ for $i<n$ and $\alpha_n=e_{n-1}+e_n$. Choose a compatible hyperbolic basis so that the chosen ruling contains the residual plane labelled by $-n,-(n-1)$. In signed one-line notation, the minimal parabolic representative is
\begin{equation}\label{eq:D-word}
 w=(p,p+2,p+3,\ldots,r,-n,-(n-1),1,2,\ldots,p-1,p+1),
\end{equation}
omitting empty strings. The first positive block contributes $(k-2)p-1$ inversions with the final positive block; the two negative entries contribute $2p$ negative images of difference roots and one negative image of a sum root. All other contributions vanish. Thus $\ell(w)=kp$, and its unique right descent is $k$. The corresponding fixed $k$-plane belongs to the incidence image, so its Schubert closure, of dimension $kp$, is $X$. The geometric proof of Lemma~\ref{D:image} does not require this label.

\subsection{The canonical obstruction along the contracted locus}
Use additive line-bundle notation and set $H=\det\mathcal L^*$ and $\zeta=\OO_{B_0}(1)$. The relative Grassmann tangent bundle is $\mathcal L^*\otimes\mathcal Q$. Since $\det\mathcal Q=H-2\zeta$, taking its determinant and then adding $K_{B_0}=-2\zeta$ gives $K_Z=-(r+2)H+(2k-2)\zeta$. The determinant section in \eqref{eq:D-incidence} has line bundle $\det\mathcal Q\otimes(\det P)^*=H-2\zeta$. Adjunction therefore yields
\[
 K_I=-(r+1)H+(2k-4)\zeta.
\]
Here adjunction applies to a Cartier hypersurface in a smooth variety, even at the singular points of the hypersurface \cite[Chapter~II, Section~8]{Har77}\cite[Chapter~3]{BH93}.

\begin{proposition}\label{D:nonQ}
The ambient variety $X$ is non-$\Q$-Gorenstein at every scheme point of $T$ and is Gorenstein on $X\setminus T$.
\end{proposition}
\begin{proof}
The complement is isomorphic to an open subset of the Gorenstein hypersurface $I$. For $x\in T$, its fibre curve is $C_x\simeq\PP^1_{\kappa(x)}$, and $H$ and $\zeta$ have degrees zero and one on it. Thus
\[
 \deg(K_I|_{C_x})=2k-4>0.
\]

Suppose $mK_X$ were Cartier near $x$ for some $m>0$. Choose canonical Weil divisors using the same rational top form on the common function field of $I$ and $X$. The map is an isomorphism outside subsets of codimension at least two in both varieties and has no exceptional prime divisor. Hence $mK_I=q^*(mK_X)$ as Weil divisors over this neighbourhood. Both sides are Cartier on the normal variety $I$, so their line bundles agree. The right side has degree zero on $C_x$, while the left side has degree $m(2k-4)>0$, a contradiction. This argument works over $\kappa(x)$ and therefore covers nonclosed points as well.
\end{proof}

\begin{proof}[Proof of Theorem~\ref{thm:mainD}]
Set $S=q(I_2)$. It is irreducible. The exceptional locus has codimension at least five, whereas $I_2$ has codimension three. Its generic point is therefore in the isomorphism locus, so $S$ has codimension three and generic transverse node. The inequality $\dim(L\cap P)\ge2$ for $L\in T$ gives $T\subset S$. Outside $S$, the map identifies $X$ with a smooth open subset of $I$. The dense subset $S\setminus T$ is singular in $X$, and closedness of the singular locus then gives $\Sing X=S$.

Proposition~\ref{D:nonQ} gives $\NGor(X)=\NonQGor(X)=T$, with $\codim_XT=2k$. The generic transverse node along $S$ is Gorenstein, so $\Bbad(X)=\varnothing$.

The inverse image of $X$ under $G/B\to G/P$ is a complete-flag Schubert variety. This projection is smooth and locally trivial with geometrically irreducible flag variety fibre. The singular, Gorenstein, and $\Q$-Gorenstein loci pull back, and codimensions are unchanged. For the last assertion, the canonical sheaf on a local product is the pullback of the base canonical sheaf tensored with an invertible relative canonical sheaf. Reflexive powers commute with the flat extension, and invertibility is reflected by faithful flatness. Complete-flag counterexamples follow.
\end{proof}

\subsection{An independent toric model at a general point of the locus}
At $L_0\in T$ in the open position $\dim(L_0\cap P)=2$, choose decompositions $R=L_0\oplus W$ and $P=(L_0\cap P)\oplus W$. A Grassmann graph chart writes $L=\operatorname{graph}(\alpha,M)$, where $\alpha:L_0\to W$ and $M:L_0\to H_0$. The $k(r-k)$ coordinates of $\alpha$ are unconstrained. The ruling condition writes $M=x\otimes B$, with $x\in\C^2$ and $B$ a $2\times k$ matrix. Choosing the first two basis vectors in $L_0\cap P$, the incidence condition is $\det B_{\{1,2\}}=0$. Thus the actual chart is
\begin{equation}\label{eq:D-toric-chart}
 \A^{k(r-k)}\times C_k,\qquad
 C_k=\Spec\!\left(\left(\C[x_1,x_2,B_{bi}]/(\det B_{\{1,2\}})\right)^{\Gm}\right),
\end{equation}
where $x$ has weight $+1$ and $B$ has weight $-1$.

To check the scheme identity, the ring before invariants is a rank-four quadric ring times a polynomial ring, hence is normal and integral. Its invariant ring is likewise normal and integral. Pairing positive and negative factors shows that it is generated by the entries $x_aB_{bi}$, so it defines a closed reduced subscheme in the graph-coordinate affine space. Its closed-point support is the incidence image: a nonzero residual map recovers the unique ruling line, and the zero map is represented as well. Both this subscheme and the chart of $X$ are reduced, so equality of their closed-point supports over $\C$ gives equality of defining ideals. The scheme-theoretic identification does not require the equations $x_ax_c\det B_{\{1,2\}}$ to generate the entire defining ideal.

Equivalently, $C_k$ is the full $\OO(1,1)$ section cone over $\PP^1\times Q_{2k-2}$, where $Q_{2k-2}\subset\PP^{2k-1}$ is the rank-four quadric $\det B_{\{1,2\}}=0$. Its degree-$d$ dimension is
\[
 (d+1)\left[\binom{d+2k-1}{2k-1}-\binom{d+2k-3}{2k-1}\right],
\]
with binomial coefficients zero when appropriate. Applying $1+t\frac{d}{dt}$ to $(1+t)/(1-t)^{2k-1}$ gives
\[
 H_{C_k}(t)=\frac{1+2kt+(2k-3)t^2}{(1-t)^{2k}}.
\]
This is an identity for every $k$, not an extrapolation from examples.

Factor the first two columns of $B$ as $v\otimes w$, and call the other $2k-4$ entries $d_j$. The toric monoid has nonnegative coordinates $(u_1,u_2,v_1,v_2,w_1,w_2,d_1,\ldots,d_{2k-4})$, subject to
\[
 V=W,\qquad U=V+D,
\]
where capital letters denote coordinate sums. It is generated by the eight monomials $u_av_bw_c$ and the $4k-8$ monomials $u_ad_j$. Pairing factors in these balance equations proves that this is exactly the nonnegative cone intersected with the indicated lattice, so the monoid is saturated. This lattice is also the group generated by the monoid: the solution $u_1=u_2=k-1$ with all $v_i,w_i,d_j$ equal to one is strictly positive, and adding a sufficiently large multiple of it expresses any lattice vector as a difference of two nonnegative solutions. Each coordinate inequality is a primitive facet: its zero face contains a point with all other coordinates positive, and one of the generators has coordinate value one. The toric canonical divisor is $K=-\sum_\rho D_\rho$ \cite[Theorem~8.2.3]{CLS11}.

If $mK_{C_k}$ were Cartier at the vertex, its divisorial module would have a one-dimensional quotient by the vertex ideal. That quotient is torus-stable. Lift a weight vector using diagonalizability, then apply the positive-degree graded Nakayama lemma to obtain a generator of the whole module. In the full torus grading each character space has dimension at most one, so this generator is a monomial. Ordinary cone homogeneity alone would not justify this last assertion. At every invariant prime divisor, the valuation threshold of $\OO(mK)$ is $m$, so the primitive-facet value of the generator is exactly $m$. All these divisors pass through the vertex. Substitution into the second balance equation forces $2m=(2k-2)m$, impossible for $k\ge3$ and $m>0$.

This independently recovers the infinite-order canonical obstruction. Under the lattice map with coordinates $(V-W,U-V-D)$, the canonical vector maps to $(0,2k-4)$, in agreement with the adjunction calculation. For $(r,k)=(3,3)$, the local six-dimensional cone has Hilbert numerator $(1,6,3)$. For $(r,k)=(5,4)$, the transverse eight-dimensional cone along $T\simeq\PP^4$ has numerator $(1,8,5)$.

\section{Finite certificates and reproducibility}\label{sec:reproducibility}

The type~A proof and the type~D incidence construction have no finite enumeration as a mathematical premise. The arbitrary-type torsion theorem uses the bounded inputs of Proposition~\ref{tor:certificates}, after the geometric Levi reduction. The supplied programs reconstruct these inputs from the simple reflections using the Python standard library and exact integer arithmetic. The optional type~A diagnostics are finite corroboration of the geometric proof, not a substitute for it or a formal proof certification.

\subsection{The finite inputs}
For each element $w$ of a bounded Weyl group, choose one reduced word $Q_w$. One word suffices: every local ring needs one suitable regular proper birational source, not a degree bound for every reduced expression. Its full vertex-divisor column set is
\[
 E(Q_w,u)=\{j:\Dem(Q_w\text{ with position }j\text{ deleted})=u\}.
\]
Empty-column cases are covered by Theorem~\ref{tor:degrees} and Lemmas~\ref{tor:fibre}--\ref{tor:columns}.

Breadth-first enumeration in the simple-root reflection matrices constructs the entire group and a first-arrival reduced word for every element. Closure under all simple generators and the expected group order are checked. Every single deletion is evaluated with the Demazure rule and every witness mask with ordinary multiplication. The program retains the literal integer simple-root coordinates of the boundary-character differences from Proposition~\ref{tor:mask-degrees}.

For a pair of gap $d$, a mask ending at $u$ can omit at most $d$ letters: the length of its selected word is at least $\ell(u)$. The parity of the number omitted equals the parity of $d$. This proves the finite mask bound used in the reconstruction; it is not a heuristic truncation. Each witness minor uses every exceptional column and distinct legal row indices. Extended Euclid supplies explicit integral B\'ezout coefficients for the retained maximal minors, and the program checks the resulting identities.

\begin{table}[htbp]
\centering
\small
\begin{tabular}{lrrl}
\toprule
Type and gap & Group order & Centre pairs & Result\\
\midrule
$A_3$, gap $3$ & $24$ & $1$ & primitive\\
$B_3$, gap $3$ & $48$ & $12$ & all primitive\\
$C_3$, gap $3$ & $48$ & $12$ & eleven primitive; one gcd $2$\\
$C_2$, gap $2$ & $8$ & $3$ & two primitive; one gcd $2$\\
$F_4$, gap $4$ & $1152$ & $567$ & all primitive\\
$G_2$, gap $3$ & $12$ & $0$ & empty-column case\\
\bottomrule
\end{tabular}
\caption{The complete bounded inputs reconstructed by the supplied verifier. Pair counts concern nonempty exact-column sets at the indicated gap. The word argument in Proposition~\ref{tor:exceptional} proves the absence of $G_2$ columns at every gap at least three.}
\label{tab:certificate-summary}
\end{table}

For both $B_3$ and $C_3$ there are $199$ exact-centre pairs in total: $138$ of gap one, $49$ of gap two, and $12$ of gap three. Thus the exceptional pairs of gap at least two number $61$. The distinction is relevant to coverage even though the theorem uses only the specified dimensional inputs. The final simple roots are $e_3$ and $2e_3$, respectively; sharing an abstract Weyl group does not permit sharing their character matrices.

In $F_4$, the complete $567$ gap-four centre pairs comprise $497$ one-column, $49$ two-column, and $21$ three-column cases. The program enumerates all centres directly from the simple reflections and reconstructs every endpoint, character row, exact-column set, determinant, and B\'ezout identity. The resulting witnesses cover all $567$ cases.

A second implementation acts on explicit Euclidean root sets by permutations. Without importing the Cartan-matrix implementation, it checks word lengths by positive-root inversion counts, recomputes the complete deletion-centre sets and all stored character witnesses, and evaluates the witness determinants by the Leibniz formula. This supplies an independent check on the root-length conventions and the certificate calculations.

Containment of the witness rows in the actual integral curve-degree lattice, proved in Proposition~\ref{tor:mask-degrees}, is sufficient. Generation of the whole degree lattice is unnecessary. A primitive full-column submatrix forces the ramification congruences. For a general full-rank matrix, the maximal-minor gcd is the product of the nonzero Smith invariants, not necessarily the largest invariant. In these certificates the primitive cases have gcd one; the nonprimitive cases have a single column with gcd two. The verifier uses precisely those justified conclusions.

\subsection{Files and reproduction}
The programs, exact certificates, and execution records are available among the ancillary files of this arXiv version. Downloading the complete article source preserves their directory structure and their correspondence with the manuscript. The \texttt{anc/} directory contains \texttt{README.md}, the certificates in \texttt{generated/}, and the following entry points:
\begin{center}
\small
\begin{tabular}{@{}ll@{}}
\toprule
File in \texttt{anc/} & Purpose\\
\midrule
\texttt{verify\_all.py} & Run all required checks\\
\texttt{verify\_torsion\_certificates.py} & Reconstruct the exact certificates\\
\texttt{verify\_root\_permutations.py} & Check with explicit root permutations\\
\bottomrule
\end{tabular}
\end{center}

Python~3.10 or later and its standard library suffice. From the source-package root, use the first command for the required checks and the second to include the optional type~A diagnostics:
\begin{quote}
\small
\begin{verbatim}
python anc/verify_all.py
python anc/verify_all.py --extended
\end{verbatim}
\end{quote}
The default run reconstructs all six rows of Table~\ref{tab:certificate-summary}, compares the supplied JSON with the reconstruction, independently checks the root-permutation witnesses, and checks the type~D sample identities. Any failed check gives a nonzero exit status. Success ends with \mbox{\texttt{ALL REQUESTED CHECKS PASSED.}} The certificate manifest separately checks file identities.

Measured in a fresh standard-library virtual environment on Linux with CPython~3.12.14 and an AMD EPYC~9V74 processor, the default and extended runs took $2.3$ and $33.3$ seconds, respectively. These single-run wall-clock measurements, together with the environment details, individual timings, and complete outputs, are recorded in \texttt{anc/verification\_results/}.

The extended run checks $417$ filtered-maximum instances in $S_6$, $4{,}998$ in $S_7$, and, in $S_6$, $347{,}542$ horizontal-incidence comparisons against an independent bipartite matching algorithm, $29{,}968$ greedy lifts, and $4{,}096$ complete horizontal-prefix witnesses. These bounded checks do not establish the arbitrary-rank theorem; the coordinate and divisor proofs in Section~\ref{sec:type_a} do that.

All arbitrary-rank and field-extension quantifiers come from the geometric arguments. In particular, the type~C scalar-Veronese identification uses the actual Lagrangian graph model and the equivariant induction. It is not inferred from a Smith invariant equal to two. Likewise, the $D$-family dimension, canonical degree, and Hilbert-series identities are proved algebraically for every parameter in Section~\ref{sec:type-d}; the finite sample checks serve only to detect transcription errors.

\section*{Acknowledgments}
I thank Henning Ulfarsson for sharing the Permutation Pattern Avoidance Library (PermPAL; \url{https://permpal.com}), a resource for exploring permutation pattern classes.

The proofs in this paper were generated with the assistance of GPT.
\bibliographystyle{alpha}
\clearpage
\newcommand{\etalchar}[1]{$^{#1}$}

\end{document}